\RequirePackage[l2tabu, orthodox]{nag}
\PassOptionsToPackage{table,xcdraw}{xcolor}

\newif\ifllncs
\newif\ifspringer
\newif\ifnonspringer

\newif\ifanonymous
\anonymousfalse

\llncsfalse
\springerfalse
\nonspringertrue

\makeatletter
\providecommand*{\input@path}{}
\g@addto@macro\input@path{{llncs/}{include/}{fig/}{springer/}}
\makeatother

\ifllncs
    \documentclass[runningheads]{llncs}
    \usepackage[misc]{ifsym} %
    \ifanonymous
        \author{}
        \institute{}
    \else
        \author{
            Boyang Han\inst{1}\orcidID{0009-0009-1527-2336} \and\\[0pt]
            Aleksandr M. Kazachkov\inst{1}\textsuperscript{(\Letter)}\orcidID{0000-0002-4949-9565}
        }
        \authorrunning{B.~Han and A.~M.~Kazachkov}
        
        \institute{
            Center for Applied Optimization,
            Department of Industrial \& Systems Engineering,
            University of Florida, Gainesville, FL 32611, USA\\
            \email{bhan1@ufl.edu, akazachkov@ufl.edu}
        }
    \fi
\fi

\ifspringer
    \RequirePackage{amsmath}
    \documentclass[smallextended,envcountsame,nospthms]{springer/svjour3}
    \journalname{Mathematical Programming}
    \smartqed  %
    \usepackage{marvosym} %
    \usepackage{etoolbox}
    \makeatletter
    \patchcmd{\acknowledgement}{\addvspace{17pt}}{\addvspace{8pt}}{}{} %
    \makeatother
    
    \author{Boyang Han \and Aleksandr M. Kazachkov}

    \institute{%
           B. Han (ORCID 0009-0009-1527-2336)
           \and
           A. M. Kazachkov (ORCID 0000-0002-4949-9565)
             \at
            Center for Applied Optimization,
            Department of Industrial \& Systems Engineering,
            University of Florida, Gainesville, FL 32611, USA\\
            \email{bhan1@ufl.edu, \Letter{} akazachkov@ufl.edu}
    }
    \date{} %
    \authorrunning{B.~Han and A.~M.~Kazachkov}
\fi

\ifnonspringer
    \documentclass[11pt]{article}
    \usepackage[margin=1in]{geometry}
    \author{Boyang Han \and Aleksandr M. Kazachkov}
    \date{}
\fi

\title{The Strength of Root Cuts in an Extended Abstract Branch-and-Cut Model%
        \thanks{Extension of paper from The 27th Conference on Integer Programming and Combinatorial Optimization~\cite{HanKaz26_abstractbc-ipco}.}
}

\usepackage[T1]{fontenc}

\usepackage{amsfonts,amssymb}
\ifllncs
\else
    \RemoveFromHook{package/amsthm/after}[firstaid/aliascounter] %
    \usepackage{amsthm}
\fi
\usepackage{silence}
\usepackage{thmtools}
\usepackage{thm-restate}
\usepackage{float}
\usepackage{graphicx}
\usepackage{listings}
\usepackage{appendix}
\usepackage{textcomp}
\usepackage{tikz}
\usepackage{mathtools}

\newcommand{\w}{\mathrm{w}}
\newcommand{\W}{\tilde{\mathrm{w}}}
\newcommand{\septime}{\mathrm{s}}        %
\newcommand{\sepinf}{\septime_\infty}
\newcommand{\sepmax}{S_{\max}}
\newcommand{\nleaf}{\lambda}            %
\usepackage{longtable,outlines}
\ifspringer\else\usepackage{caption}\fi %
\usepackage{booktabs}

\usepackage{subfigure}
\usepackage{color,enumerate,mathrsfs,algorithm,algorithmicx,algpseudocode}
\usepackage{bbm} %
\usepackage{xspace}
\usepackage[dvipsnames]{xcolor}
\usepackage[shortlabels,inline]{enumitem}
\setlist{nolistsep}

\usepackage{etoolbox}
\patchcmd{\thebibliography}{\chapter*}{\section*}{}{}

\usepackage[numbers,compress]{natbib}

\usepackage[colorlinks=true,allcolors=blue]{hyperref}
\usepackage[nameinlink,capitalize,noabbrev]{cleveref}
\definecolor{mycolor}{RGB}{200,0,0}

\usepackage{tikz}
\usetikzlibrary{calc}
\usepackage{ifthen}
\usetikzlibrary{arrows, arrows.meta, shapes, backgrounds}
\usetikzlibrary{decorations.pathmorphing, patterns} %
\usetikzlibrary{graphs, positioning, fit, quotes}
\usepackage{forest}
\tikzset{itria/.style={
        draw,shape border uses incircle,
        isosceles triangle,shape border rotate=90,yshift=-17pt}}
\tikzset{empty/.style={draw=none,fill=none}}
\tikzset{cut/.style={regular polygon, regular polygon sides=4,minimum size=2.1213mm}}
\tikzset{branch/.style={circle,minimum size=1.5mm}}
\tikzset{big cut/.style={regular polygon, regular polygon sides=4,minimum size=7.0711mm}}
\tikzset{big branch/.style={circle,minimum size=5mm}}
\forestset{small node/.style={inner sep=0pt, branch, fill=black, draw, font=\small},
    l=10pt,->}
\forestset{big node/.style={inner sep=0pt, big branch, draw, font=\small},
    l=10pt,->}
\forestset{
    join aunts/.style={
    before drawing tree={
      tempkeylista'=,
      for nodewalk={fake=u, siblings}{tempkeylista/.option=name},
      join list/.register=tempkeylista,
      tikz+/.process={ OOw2 {join list} {fork sep} { \draw  [thick, rounded corners, Stealth-]\foreach \i in {##1} { (.child anchor) -- ++(0,##2) -| (\i.parent anchor) }; } },
    }
  }
}
\tikzset{
  subtreewt/.append  style={align=center,
    label={[fill=none,font=\small,yshift=-0.25cm]0:{\hspace{1.5pt} subtree \heavycutweight = #1}}},
}
\tikzset{
  subtreewtdeep/.append  style={align=center,
    label={[fill=none,font=\small,yshift=-0.6cm]0:{\hspace{2.5pt} subtree \heavycutweight = #1}}},
}
\tikzset{
  subtreewtdeep_new/.append style={
    align=center,
    label={[fill=none, font=\small, yshift=-0.6cm, xshift = 0.3cm, align=center]0:{subtree \heavycutweight \\ = #1}}
  },
}

\tikzset{
  lvlwt/.append  style={align=center,
    label={[fill=none,font=\small]0:{\hspace{1.5pt} layer \heavycutweight = #1}}},
}
\tikzset{
  plvlwt/.append  style={align=center,
    label={[fill=none,font=\small]0:{\hspace{1.5pt} PLT = #1}}},
}
\tikzset{
  nodewt/.style={align=center,
    label={[fill=none,font=\small]0:{\hspace{1.5pt} \heavycutweight = #1}}},
}

\tikzset{subtree/.style={draw, isosceles triangle, shape border rotate=90,
    isosceles triangle apex angle=50, anchor=apex, inner sep=0.5pt,
    font=\scriptsize, minimum height=#1}}
\forestset{
  tri/.style 2 args={tikz+={\node[subtree=#2] at (.south) {\raisebox{-1.2ex}{$#1$}};}},
  cutpath/.style={edge={dotted,thick}},
}

\makeatletter
\DeclareRobustCommand{\rvdots}{%
  \vbox{
    \baselineskip4\p@\lineskiplimit\z@
    \kern-\p@
    \hbox{.}\hbox{.}\hbox{.}
  }}
\makeatother
\makeatletter
\def\Ddots{\mathinner{\mkern1mu\raise\p@
\vbox{\kern7\p@\hbox{.}}\mkern2mu
\raise4\p@\hbox{.}\mkern2mu\raise7\p@\hbox{.}\mkern1mu}}
\makeatother
\forestset{dot node/.style={big node,font={$\rvdots$},draw=none,fill=none}}

\newcommand{\dotsangle}{26}
\newcommand{\dotsxshift}{.2ex}
\newcommand{\dotsyshift}{1ex}
\newcommand{\Ldots}{\hspace{\dotsxshift}%
    \raisebox{\dotsyshift}{\rotatebox{\dotsangle}{$\Ddots$}}}

\forestset{ldiag dot node/.style={big node,font={$\Ldots$},draw=none,fill=none}}
\forestset{rdiag dot node/.style={big node,font={$\ddots$},draw=none,fill=none}}

\renewcommand{\tilde}{\widetilde}
\renewcommand{\hat}{\widehat}

\newcommand*{\defeq}{\mathrel{\vcenter{\baselineskip0.5ex \lineskiplimit0pt
                     \hbox{\footnotesize.}\hbox{\footnotesize.}}}%
                     =}

\colorlet{rosso}{red!80!blue}

\usepackage[normalem]{ulem}

\newcommand{\cardinality}[1]{\lvert #1 \rvert}

\ifnonspringer
    \newtheorem{theorem}{Theorem}
    \newtheorem*{theorem*}{Theorem}
    
    \newtheorem{corollary}[theorem]{Corollary}
    \newtheorem{lemma}[theorem]{Lemma}
    \newtheorem{proposition}[theorem]{Proposition}
    \newtheorem{claim}[theorem]{Claim}
    \newtheorem{remark}[theorem]{Remark}
    
    \newtheorem{conjecture}[theorem]{Conjecture}

    \theoremstyle{remark}
    
    \newtheorem*{subclaim*}{Claim}

    \theoremstyle{definition}
    
    \newtheorem{definition}[theorem]{Definition}
    \newtheorem*{definition*}{Definition}
    \newtheorem{example}{Example}
    \newtheorem*{example*}{Example}
    
\fi

\ifspringer
    \newtheorem{example}{Example}
    \newtheorem{lemma}[example]{Lemma}
    \newtheorem*{lemma*}{Lemma}
    \newtheorem{proposition}[example]{Proposition}
    \newtheorem*{proposition*}{Proposition}
    \newtheorem{theorem}[example]{Theorem}
    \newtheorem*{theorem*}{Theorem}
    \newtheorem{corollary}[example]{Corollary}
    
    \newtheorem{definition}[example]{Definition}

    \newtheorem*{claim*}{Claim}
    \theoremstyle{remark}
    
\fi

\allowdisplaybreaks

\ifspringer
    \Crefname{figure}{Fig.\!}{Figs.\!}
    \Crefname{section}{Sect.\!}{Sects.\!}
    \Crefname{definition}{Def.\!}{Defs.\!}
\fi
\crefname{enumi}{Statement}{Statements}
\crefname{line}{Step}{Steps}
\crefname{equation}{}{}

\newcommand{\suchthat}{:}
\newcommand{\R}{\mathbb{R}}

\newcommand{\Z}{\mathbb{Z}}
\newcommand{\reals}{\R}
\newcommand{\integers}{\Z}
\newcommand{\nonnegreals}{\reals_{\scriptscriptstyle \ge 0}}
\newcommand{\nonnegints}{\integers_{\scriptscriptstyle \ge 0}}
\newcommand{\posreals}{\reals_{\scriptscriptstyle > 0}}
\newcommand{\posints}{\integers_{\scriptscriptstyle > 0}}

\newcommand{\bb}{B\&B\xspace}
\newcommand{\bc}{B\&C\xspace}
\newcommand{\cp}{CP\xspace}
\renewcommand{\bar}{\overline}
\renewcommand{\hat}{\widehat}
\renewcommand{\tilde}{\widetilde}

\newcommand{\treeweight}{\operatorname{\tau}}
\newcommand{\treetime}{\treeweight}

\newcommand{\heavycutweight}{time\xspace}

\newcommand{\lbd}{\ell}
\newcommand{\rbd}{r}
\newcommand{\cutbd}{c}
\newcommand{\cbd}{\cutbd}
\newcommand{\cutfunset}{\mathcal{F}_{\scriptscriptstyle\searrow}}
\newcommand{\ncutfams}{m}

\newcommand{\targetbound}{Z}
\newcommand{\targetbd}{\targetbound}

\newcommand{\nodes}{\mathcal{V}}

\newcommand{\historyvec}{\vec{h}}

\newcommand{\nodetypes}{\mathcal{I}}

\NewDocumentCommand{\PureBB}{o m}{%
  T\IfNoValueTF{#1}{(#2)}{(#1,#2)}%
}

\newlength{\subplotwidth}
\usepackage{color}

\begin{document}

\maketitle              %
\begin{abstract}
The branch-and-cut paradigm has led to effective solvers for a wide range of mathematical optimization problems. Despite significant theoretical work on pure branching and pure cutting, less is known about their interaction. We study an extended abstract branch-and-cut model to investigate tradeoffs between branching and cutting, including the time to generate cuts. We prove that, under symmetric branching, some minimal-time tree has all cuts at the root, before the first branching node, and that there is an asymmetric setting where some minimal-time tree has no internal cuts, i.e., cuts both preceded and followed by branching. For other settings, although a root-cut tree may be suboptimal and differ substantially in structure from a minimal-time tree, it nevertheless provides a reliable approximation whose tree time does not deviate arbitrarily from the optimum. For constant cut strength and bounded time functions, we show that the best root-cut tree is asymptotically optimal as the target bound grows. In addition, under a single-cut-per-path restriction, we characterize when some minimal-time tree has an internal cut. Lastly, simulations indicate that the best root-cut tree is much closer to a minimal-time tree than our worst-case approximation guarantee suggests.
   \ifnonspringer\else
   \fi
   \ifspringer
   \fi
\end{abstract}
\section{Introduction}
Branch-and-cut (\bc), a hybrid version of the branch-and-bound (\bb)~\cite{LanDoi60} and cutting plane (\cp) methods~\cite{Gomory58_fractional-cut,Gomory60_solving-linear-programming-problems-in-integers,Gomory63_algorithm-for-integer-solutions-to-linear-programs}, is a powerful tree-based framework 
that is a key ingredient of mathematical optimization solvers.
Each node of the tree corresponds to a relaxation of the original optimization problem that is refined by either valid inequalities (\emph{cuts}) or splitting the search space (\emph{branching}).
Solvers determine \emph{where} in the tree to add rounds of cuts from empirical guidelines and typically focus on the \emph{root node} of the tree.
In this paper, we provide complementary theory on cut placement with provable guarantees on the computational cost of the resulting tree, within an extension of the \emph{abstract \bc model} by
Kazachkov, Le~Bodic, and Sankaranarayanan
\cite{KazLeBSan22_abstractbc-ipco,KazLeBSan24_abstractbc-mpb}.

Our model assumes that key solution parameters are known:
    the bound that needs to be proved (e.g., an optimal solution is given, but we need to certify optimality via the tree),
    the evolution of the bound over iterated branching and cutting operations,
    and an estimated time to solve each node of the tree.
Our perspective (when such information may be available) is a \emph{post-optimality analysis} of the quality of the best-possible \bc tree for a problem, 
which can provide beneficial insights when solving the same problem repeatedly.

In the original abstract \bc model~\cite{KazLeBSan24_abstractbc-mpb}, branching yields fixed (constant) bound improvements, while cuts may have diminishing strength across rounds or increase node processing time. Under symmetric branching and constant cut strength, they prove that a \emph{root-cut tree}, where all cut nodes are before the first branch node, is sufficient:
there always exists a minimal-time root-cut tree,
for any function defining node solving time~\cite[Theorem~23]{KazLeBSan24_abstractbc-mpb}.

We generalize this model to accommodate broader branching and cutting scenarios 
and show that a root-cut tree solving in minimal (or nearly minimal) time continues to exist.
Even when the root-cut tree is suboptimal, 
i.e., when any minimal-time tree requires \emph{local} cuts (deeper in the tree),
we derive an upper bound 
on how much its processing time can deviate from the minimal-time tree,
potentially offering guidance (or validation) for the root-cut dominant approach of solvers.
Our contributions are summarized next.

\begin{itemize}[label=$\triangleright$,leftmargin=*,after=\vspace{-1ex}]
    \item \Cref{sec:root}: We prove that, under symmetric branching, 
    a minimal-time root-cut tree always exists,
    even when simultaneously considering multiple cut and branching families, or node-dependent cut generation times. %
  \item \Cref{sec:sparse-then-dense}: We prescribe how to best add a sequence of weak and strong cuts, when the weak cuts cause a relatively smaller increase in node solving time. 
  \item \Cref{sec:124}: 
    We identify an asymmetric branching setting where some minimal-time tree has no internal cuts, i.e., all cuts are at the root or at the leaves.
  \item \Cref{sec:notbad}: We bound the approximation error of the best root-cut tree with respect to the minimal-time tree. %
  \item \Cref{sec:asymptotic}: 
  We analyze bounded time functions, for which the best root-cut tree is within a constant factor of a minimal-time tree,
  and for constant cuts, it is asymptotically optimal.
  \item \Cref{sec:complex}: We show by example that the minimal-time tree can be structurally complex for asymmetric settings: a local cut may appear at any depth.
  \item \Cref{sec:1cut}: We characterize conditions that allow the construction of minimal-time trees featuring local cuts under a single-cut-per-path assumption.
  \item \Cref{sec:separation}: We study how cut generation time impacts number of cut nodes in a minimal-time tree, proving monotonicity for constant separation times.
  \item \Cref{sec:simulation}: We empirically assess the quality of the root-cut approximation. %
\end{itemize}

\paragraph{Related literature.}
Recent theory has been developed to analyze the strength of the pure \bb{} process \cite{DeyDubMol22_lower-bounds-on-size-of-general-branch-and-bound-trees_mp,DadTiw20_on-the-complexity-of-branching-proofs_ccc,DeyDubMolSha24_analysis-of-full-strong-branching_mp,LeBNem17}, pure \cp \cite{DeyMol18}, and comparing \bc, \bb, and \cp~\cite{BasConSumJia22_complexity-of-bb-and-cuts-II,BasConSumJia23_complexity-of-bb-and-cuts-I}. 
Basu, Conforti, Di~Summa, and Jiang~\cite{BasConSumJia22_complexity-of-bb-and-cuts-II} show that 
\bc can be strictly superior to \bb and \cp. 
In contrast, and counterintuitively, 
Shah, Dey, and Molinaro~\cite{ShaDeyMol26_nonmonotonicity-branching-rules_ijoc} prove that a single cut can make \bb exponentially worse. 
Related to our work, Jiang~\cite{Jiang26_sharp-logarithmic-thresholds_arxiv} identifies---within our model---the fraction of the target bound proved by cut nodes in a minimal-size \bc{} tree, as a function of how cumulative cut strength grows.
Other papers adopt a learning-based perspective:
Berthold, Francobaldi, and Hendel~\cite{BerFraHen25_learning-to-use-local-cuts} propose methods for learning to use local cuts effectively within \bc{};
Cheng and Basu~\cite{CheBas25_generalization-guarantees-for-score-based-b-and-c-policies_neurips} establish generalization guarantees for robustness of learned \bc policies
across problem instances.
However, in all this research, the explicit tradeoffs and interplay between branching and cutting decisions remain poorly understood.
Our high-level motivation is to systematically study the interaction between branching and cutting within a \bc framework in a stylized (but already nontrivial to analyze) model, through which we can progress towards theoretical guidance for cut placement and selection in solvers. %

\section{Preliminaries: Notation and Definitions}\label{notion}

Abstractly, \bc involves solving a sequence of tightened subproblems in a tree structure. For example, in mixed-integer optimization, the \bc process starts from the \emph{root node}, a convex relaxation of the original problem. The optimal value at the root node is a dual bound for the original problem. Then, the problem formulation (and dual bound) is strengthened by either: (1) partitioning the feasible space and creating two subproblems, or (2) adding a round of valid cuts that are violated by points in the current relaxation but satisfied by all feasible mixed-integer solutions. 
The \bc tree is built by recursively selecting between these two operations until the tree \emph{proves a (dual) bound} matching the primal optimal value (we assume this is given as $\targetbound$).

We largely adopt the notions modeling a \bc tree and its processing time from Kazachkov, Le~Bodic, and Sankaranarayanan \cite[Definitions 1–5]{KazLeBSan24_abstractbc-mpb}.
The authors define a \emph{time function}~\cite[Definition~2]{KazLeBSan24_abstractbc-mpb} to model that cuts make formulations more complex, thus increasing the processing time of subsequent subproblems.
In practice,
    solvers can choose among \emph{different families of cuts} that have heterogeneous properties,
    and \emph{cut generation time} also affects cut selection.
We introduce two extended models of solving time:
a \emph{simple} time function (\Cref{def:simple-time-function}) in which solving time of a node in the tree only depends on the number of cuts previously generated from each family, 
and a \emph{generalized} time function (\Cref{def:timefun-generalized}) that also accounts for the time to compute cuts.
We will assume that no cut family strictly dominates another (i.e., when a cut family yields a better bound improvement, it causes slower relaxations). 

Let $\nodes_T$ denote the nodes of a \bc tree $T$.
In our model, each node in the tree has a type $i \in \nodetypes \defeq \{0,1,\ldots,\ncutfams\}$,
where type $0$ represents a branch node, and type $i \in [\ncutfams] \defeq \{1,\ldots,\ncutfams\}$ represents a cut node from family $i$.
We make the simplifying assumption that the bound improvement from cut family $i \in [\ncutfams]$ only depends on the number of cuts previously added from that family.
We represent the bound improvement from a cut node from family $i$ via a \emph{(worsening) cut function} $\cbd_i \in \cutfunset$,
   where $\cutfunset$ denotes the set of nonincreasing functions $\posints \to \nonnegreals$
    and $\cbd_i(k)$ represents the bound improvement by the $k$th round of cuts from family $i$.
The bound improvement from a branch node is given by constants $\lbd$ (to the left) and $\rbd$ (to the right).
We assume that $\lbd \le \rbd$ (without loss of generality) and
that cut and branch nodes are derived from the same information (e.g., disjunction),
so $\cbd_i(k) \le \lbd \le \rbd$ for every $k \in \posints$ and cut family $i \in [\ncutfams]$.
For a node $v \in \nodes_T$, we denote by $T_v$ the subtree of $T$ rooted at $v$.
Each node $v$ also has an associated bound, denoted $z_T(v)$, defined as the total bound improvement along the path from the root of $T$ to $v$, where the root of a tree has bound 0; see \cite[Definition~1]{KazLeBSan24_abstractbc-mpb}.
A tree proves a bound $\targetbd$ if the bound at all leaves is at least $\targetbd$.

Let $\historyvec(v) \in \nonnegints^\nodetypes$ be the \emph{history vector},
where $\historyvec(v)_i$ is the number of nodes of type $i$ applied before a node $v \in \nodes_T$.
Let $e_i \in \R^\nodetypes$ denote the unit vector in direction $i$.

\begin{definition}[Simple time function]
    \label{def:simple-time-function}
    Given cut functions $\{c_i\}_{i=1}^{\ncutfams}$,
    a function $\w: \nonnegints^\nodetypes \to [1,\infty)$ is a \emph{simple time function} if it is elementwise nondecreasing, $\w(\mathbf{0}) = 1$, and if,
    for all history vectors $\historyvec \in \nonnegints^\nodetypes$, when the next cut from family $i$ is stronger than one from family $j$, i.e., $c_i(\historyvec_i+1) \ge c_j(\historyvec_j+1)$, 
    then %
    $\w(\historyvec+e_i) \ge \w(\historyvec+e_j)$.
\end{definition}

\noindent
Though \Cref{def:simple-time-function} allows dependence on the number of branch nodes, in this work,
we assume that if $\historyvec_i=\historyvec'_i$ for all $i \in \nodetypes \setminus \{0\}$, then $\w(\historyvec)=\w(\historyvec')$.
As a result, for simplicity, we omit index $0$ for a history vector input to $\w$.

\begin{definition}[Separation time function]
    \label{def:separator-time-function}
    A family of functions $\{\septime_i \suchthat i \in \nodetypes\}$ with $\septime_i: \nonnegints \to \nonnegreals$ for all $i \in \nodetypes$, is a \emph{separation time function family} if $\septime_0 \equiv 0$ and $\septime_i(\historyvec_i)$, $i \in [\ncutfams]$, represents the time to generate one round of cuts from family $i$ given $\historyvec_i$ prior cut rounds from that same family.
\end{definition}

\begin{definition}[Generalized time function]
    \label{def:timefun-generalized}
    Given a simple time function $\w$ and a separation time function family $\{\septime_i\}_{i=0}^{\ncutfams}$,
    a function $\W: \nonnegints^\nodetypes \times \nodetypes \to [1,\infty)$
    is a \emph{generalized time function} if for all history vectors $\historyvec$ and node types $i \in \nodetypes$, it holds that
        $\W(\historyvec,i) = \w(\historyvec) + \septime_i(\historyvec_i)$.  
Given a \bc tree $T$ and a node $v \in \nodes_T$ of type $i \in \nodetypes$ having a history vector $\historyvec(v)$, we will also use $\W(v) \defeq \W(\historyvec(v),i)$.
\end{definition}

Optimizing tree size, i.e., $\w \equiv 1$, is studied by Le~Bodic and Nemhauser \cite{LeBNem17} (for \bb) and Kazachkov, Le~Bodic, and Sankaranarayanan \cite{KazLeBSan24_abstractbc-mpb} (for \bc). 
However, %
a smaller \bc tree does not guarantee a better solving time. 
\Cref{ex:shorter-but-slower} 
illustrates this phenomenon.
In this work, we focus on optimizing \emph{tree time} $\treetime_{\W}(T) \defeq \sum_{v \in \nodes_T} \W(v)$ with respect to nontrivial time functions ($\W \not\equiv 1$).
The resulting optimal structure is referred to as a \emph{minimal-time tree}.
We call a tree proving bound $\targetbd$ \emph{cut-minimal} (with respect to a given time function) if it is a minimal-time tree and has the fewest cut nodes among all minimal-time trees proving $\targetbd$. 
A cut-minimal tree always exists in our model.
A \emph{pure cutting tree} has no branch nodes
and a \emph{pure branching tree} has no cut nodes.
For $z \in \reals$, let $\PureBB{z}$ denote the (unique~\cite{LeBNem17}) smallest pure branching tree proving bound $z$ for fixed $\lbd$ and $\rbd$,
in which a node is a branch node if and only if its bound is less than $z$.

Our focus is on the \emph{structure} of the tree, meaning the relative placement of cut and branch nodes.
To this end, we distinguish between \emph{root cut} nodes (applied in the tree before any branch node) and \emph{local cut} nodes (applied after the first branch node).
Define a \emph{leaf cut} node as a cut node without a branch node descendant
(i.e., part of a path ending in a \emph{leaf node}),
and an \emph{internal cut} node as a local (not root) cut that is not a leaf cut.
Define a \emph{branch-then-cut} node to be a branch node with at least one internal cut child.

\begin{example}
\label{ex:shorter-but-slower}
\Cref{fig:timevssize} shows two \bc trees that prove bound $\targetbound=9$. 
In the left panel, the tree  has 10 nodes and solves in tree time $1+ 9\w(1) = 12.7$. 
The tree in the right panel has 11 nodes, but solves more quickly, with a tree time of $7+ 4\w(1) = 12.2$,
and is in fact the unique minimal-time tree with respect to $\w$. 
In the right panel of \Cref{fig:timevssize}, the root is a branch-then-cut node because its child is an internal cut, while the node labeled 7 is a leaf cut.
\end{example}

{
    \begin{figure}
        \centering
        \subfigure{%
            \scalebox{0.8}{ %
            \begin{forest}
                for tree = {big node}
                [0,big cut 
                    [2 
                        [5 [8[11][12]][9]]
                        [6 [9][10]]
                    ]
                ]
            \end{forest}
            }
        }%
        \hspace{2cm}
        \subfigure{%
            \scalebox{0.8}{ %
            \begin{forest}
                for tree = {big node}
                [0
                    [3
                        [6 [9][10] ]
                        [7, big cut [9]]
                    ]
                    [4, big cut,
                        [6
                            [9] [10]
                        ]
                    ]
                ]
            \end{forest}
            }
        }%
        \caption{Two \bc trees proving $\targetbound = 9$, with $\cutbd = 2$, $\lbd= 3, \rbd = 4$, $\w(1) = 1.3$, and with $\w(k) = M$ for $k\ge 2$ and sufficiently large $M\in\nonnegreals$.}
        \label{fig:timevssize}
    \end{figure}
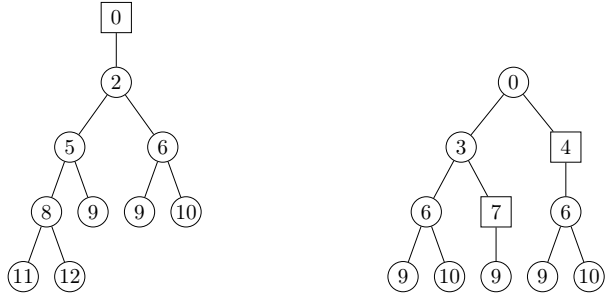
}

\section{When Root Cuts Are (Nearly) Sufficient}\label{root}

In this section, we identify settings, both symmetric (\cref{sec:root}) and asymmetric (\cref{sec:124}), in which there is a minimal-time tree without internal cuts, i.e., in which all cuts are at the root or at the leaves.
We also show (\cref{sec:sparse-then-dense}) how to best order a sequence of weak and strong cuts, when weak cuts cause a comparatively smaller increase in node solving time.

\subsection{Symmetric Branching}
\label{sec:root}

While the assumption of symmetric branching (when $\lbd=\rbd$) is restrictive in general, it is a natural theoretical starting point for our study.
Basu, Conforti, Di~Summa, and Jiang~\cite{BasConSumJia22_complexity-of-bb-and-cuts-II} (see also \cite[Example 12]{KazLeBSan24_abstractbc-mpb}) provide a concrete construction demonstrating that adding cuts can lead to an exponential improvement in the \bc tree under a symmetric branching framework.
Further, many combinatorial optimization problems, such as scheduling and set covering, have inherent algebraic or structural symmetries~\cite{Margot10_symmetry}.

We first prove \Cref{lem:nocommonprefix}, which states that no cut-minimal tree has a branch node whose children begin with the same path of cut nodes, unless both paths continue with cut nodes (from different families).%
\footnote{The extra condition is needed for our proof, but we conjecture the statement holds more generally, i.e., that no cut-minimal tree has a branch node whose two children are cut nodes of the same family.}
Ruling out this ``redundant'' substructure underlies our proofs of \Cref{thm:rootcuts,thm:124}.

\begin{lemma}\label{lem:nocommonprefix}
    Let $(\vec{\cbd},\lbd,\rbd) \in \cutfunset^\ncutfams \times \posreals \times \posreals$, and let $T$ be a cut-minimal tree proving target bound $\targetbd$ with respect to a generalized time function $\W$.
    Then $T$ contains no branch node whose two child subtrees begin with a common sequence of $k \ge 1$ cut nodes (i.e., with the same families in the same order) that, in at least one of the two subtrees, is followed by a branch or leaf node.
\end{lemma}

\begin{proof}
    Suppose, for the sake of contradiction, that
    the children of a branch node $u \in \nodes_T$ begin with identical paths of $k$ cut nodes,
    where at least one of these paths ends at a branch or leaf node.
 
 Let $T_L$ and $T_R$ be the left and right child subtrees of $u$, whose first $k$ cut nodes have the same ordered family labels. 
 Denote the cut nodes in $T_R$ by $v_1,\ldots,v_k$,
 and let $y$ be the child of $v_k$, which we may assume without loss of generality is a branch or leaf node (else, we apply an analogous argument to $T_L$).
 When $y$ is a branch node, we denote its two children subtrees as $T_{y\lbd}$ and $T_{y\rbd}$; when $y$ is a leaf, we set $\treetime_{\W}(T_{y\lbd}) = \treetime_{\W}(T_{y\rbd}) = 0$ below. The subtree remaining after the first $k$ cuts in $T_L$, denoted $T_{L}^{\mathrm{sub}}$, may begin with further cuts. Let $\historyvec(v)$ be the history vector at $v \in \nodes_T$ \emph{with respect to $T$}.
 The first $k$ cut nodes of $T_L$ have the same families and cut histories as $v_1,\ldots,v_k$, and hence the same node times.

 We first consider replacing $T_R$ in $T$ by a modified tree $T'_R$, where the path of $k$ cuts at the root of $T_R$ is shifted down a ``layer''.
 A visual aid for $T_R$ and $T'_R$ is in \Cref{fig:onecutafterbranch1}.
 The root of $T'_R$ is a branch node $u'$.
 When $y$ is a branch node, each child subtree of $u'$ starts with a path of $k$ cut nodes $v'_1,\ldots,v'_k$, followed by a copy of $T_{y\lbd}$ in the left child of $u'$ and a copy of $T_{y\rbd}$ in the right child.
 When $y$ is a leaf node, then each child of root $u'$ in $T'_R$ is a path of $k-1$ cut nodes ending with a leaf node that we will call $v'_k$ (in which case, unlike $v_k$, $v'_k$ is not a cut node and has no generation cost). Note that $T'_R$ proves a bound at least as good as $T_R$. Further, as the solving time for a node does not depend on the number of previous branch nodes, we define $\historyvec'$ as a convenient modified history vector for $T'_R$, where $\historyvec'(u') = \historyvec(u)$ and $\historyvec'(v'_i) = \historyvec(v_i)$ for $i \in [k]$. 
 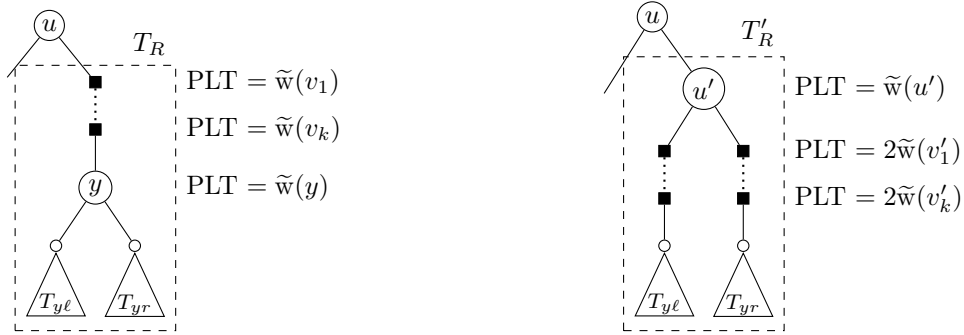
\begin{figure}[ht!]
    \centering
    \begin{minipage}[t]{0.45\textwidth}
        \centering
        \begin{forest}
    for tree = {small node, l=1pt}
    [$u$, fill=none, inner sep=1.5pt, s sep=30pt
        [, empty]
        [, cut, name=v1
            [, cut, edge={dotted, thick}, name=vk
                [$y$, fill=none, inner sep=1.5pt, s sep=25pt, name=y
                    [, fill=none, tri={T_{y\lbd}}{7mm}, name=T1]
                    [, fill=none, tri={T_{y\rbd}}{7mm}, name=T2]
                ]
            ]
        ]
    ]
    before drawing tree={
      tikz={
        \node [draw, dashed, inner xsep=4pt, inner ysep=4pt,
               fit={(v1) (T1) (T2) ([xshift=-4mm, yshift=-9mm]T1.south) ([xshift=4mm, yshift=-9mm]T2.south)}] (TR_box) {};
        \node at (TR_box.north east) [anchor=south east] {\small $T_R$};
        \node [anchor=west, font=\small] at (TR_box.east |- v1) {PLT $= \W(v_1)$};
        \node [anchor=west, font=\small] at (TR_box.east |- vk) {PLT $= \W(v_k)$};
        \node [anchor=west, font=\small] at (TR_box.east |- y)  {PLT $= \W(y)$};
      }
    }
\end{forest}
    \end{minipage}
    \hfil
    \begin{minipage}[t]{0.45\textwidth}
        \centering
        \begin{forest}
    for tree = {small node, l=1pt}
    [$u$, fill=none, inner sep=1.5pt, s sep=30pt
        [, empty]
        [$u'$, fill=none, inner sep=1.5pt, s sep=25pt, name=up
            [, cut
                [, cut, edge={dotted, thick}
                    [, fill=none, tri={T_{y\lbd}}{7mm}, name=T1]
                ]
            ]
            [, cut, name=v1
                [, cut, edge={dotted, thick}, name=vk
                    [, fill=none, tri={T_{y\rbd}}{7mm}, name=T2]
                ]
            ]
        ]
    ]
    before drawing tree={
      tikz={
        \node [draw, dashed, inner xsep=4pt, inner ysep=4pt,
               fit={(up) (T1) (T2) ([xshift=-4mm, yshift=-9mm]T1.south) ([xshift=4mm, yshift=-9mm]T2.south)}] (TR_box) {};
        \node at (TR_box.north east) [anchor=south east] {\small $T'_R$};
        \node [anchor=west, font=\small] at (TR_box.east |- up) {PLT $= \W(u')$};
        \node [anchor=west, font=\small] at (TR_box.east |- v1) {PLT $= 2\W(v'_1)$};
        \node [anchor=west, font=\small] at (TR_box.east |- vk) {PLT $= 2\W(v'_k)$};
      }
    }
\end{forest}
    \end{minipage}

    \caption{Illustration of $T_R$ (left panel) and $T'_R$ (right panel).
            PLT calculates the \emph{Partial Layer Time}, calculated only on the nodes in the visualization.
    }
    \label{fig:onecutafterbranch1}
\end{figure}

 The (sub)tree times for $T_R$ and $T'_R$ are:
 \begin{equation*}
    \begin{aligned}
        \treetime_{\W}(T_{R})
            &= \sum_{i \in [k]} \W(v_i) + \W(y) + \treetime_{\W}(T_{y\lbd}) + \treetime_{\W}(T_{y\rbd}),
    \\
        \treetime_{\W}(T_{R}')
            &= \W(u') + 2 \sum_{i \in [k]} \W(v_{i}') + \treetime_{\W}(T_{y\lbd}) + \treetime_{\W}(T_{y\rbd})
            \\
            &\le \W(u) + 2 \sum_{i \in [k]} \W(v_{i}) + \treetime_{\W}(T_{y\lbd}) + \treetime_{\W}(T_{y\rbd})
            &{\mbox{(as $\W(v_{k}') \le \W(v_{k})$).}}
    \end{aligned}
\end{equation*}
As $T$ is a minimal-time tree,
    we have that
        $\treetime_{\W}(T_{R}) \leq \treetime_{\W}(T_{R}')$,
    and we derive
        \begin{equation}
        \W(y) \leq \W(u) + \sum_{i \in [k]} \W(v_i). \label{equa:root2}
        \end{equation}

Next, let $T^{\prime\prime}$ denote a tree obtained from $T$ by shifting ``up'' a layer the path of $k$ cut nodes from both children of $u$.
Note that $T^{\prime\prime}$ proves the same bound as $T$.
Let $T_{u}$ be the subtree of $T$ rooted at $u$, 
    $T^{\mathrm{sub}}_{L}$ be the subtree of $T_{L}$ rooted at the child of the $k$th cut node, 
    and $T_{u}^{\prime\prime}$ be the tree that starts with a path of $k$ cut nodes 
        $v_{1}^{\prime\prime},\ldots,v_{k}^{\prime\prime}$ 
    followed by a branch node $u^{\prime\prime}$,
where the right child of $u^{\prime\prime}$ is the subtree $T_y$ from $T$ rooted at $y$ and the left child is $T^{\mathrm{sub}}_{L}$.
\Cref{fig:onecutafterbranch2} depicts $T_u$ and $T_{u}^{\prime\prime}$.
\begin{figure}[ht!]
    \centering
    \begin{minipage}[t]{0.45\textwidth}
        \centering
        \begin{forest}
    for tree = {small node, l=1pt}
    [$u$, fill=none, inner sep=1.5pt, s sep=30pt
        [, cut, name=L1
            [, cut, edge={dotted, thick}, name=L2
                [, fill=none, tri={T^{\mathrm{sub}}_L}{9mm}]
            ]
        ]
        [, cut, name=R1
            [, cut, edge={dotted, thick}, name=R2
                [$y$, fill=none, inner sep=1.5pt, s sep=25pt
                    [, fill=none, tri={T_{y\lbd}}{6mm}]
                    [, fill=none, tri={T_{y\rbd}}{6mm}]
                ]
            ]
        ]
    ]
    before drawing tree={
      tikz={
        \node [draw, dashed, inner xsep=10pt, inner ysep=6pt,
               fit=(L1) (L2) (R1) (R2)] (TR_box) {};
        \node at (TR_box.north east) [anchor=south west, xshift=-10pt] {\small \rlap{Sequence of $k$ cuts}};
      }
    }
\end{forest}
    \end{minipage}
    \hfil
    \begin{minipage}[t]{0.45\textwidth}
        \centering
        \begin{forest}
        for tree = {small node, l=1pt}
            [, cut, name=C1
                [, cut, edge={dotted, thick}, name=C2
                    [$u^{\prime\prime}$, fill=none, inner sep=1.5pt, s sep=25pt
                        [, fill=none, tri={T^{\mathrm{sub}}_L}{9mm}]
                        [$y^{\prime}$, fill=none, inner sep=1.5pt, s sep=22pt
                            [, fill=none, tri={T_{y\lbd}}{6mm}]
                            [, fill=none, tri={T_{y\rbd}}{6mm}]
                        ]
                    ]
                ]
            ]
        before drawing tree={
          tikz={
            \node [draw, dashed, inner xsep=10pt, inner ysep=6pt,
                   fit=(C1) (C2)] (TR_box) {};
            \node at (TR_box.east) [anchor=south west] {\small \rlap{Sequence of $k$ cuts}};
          }
        }
        \end{forest}
    \end{minipage}
    \caption{Illustration of the tree $T_u$ (left panel)
    and $T^{\prime\prime}_u$ (right panel).}
    \label{fig:onecutafterbranch2}
\end{figure}
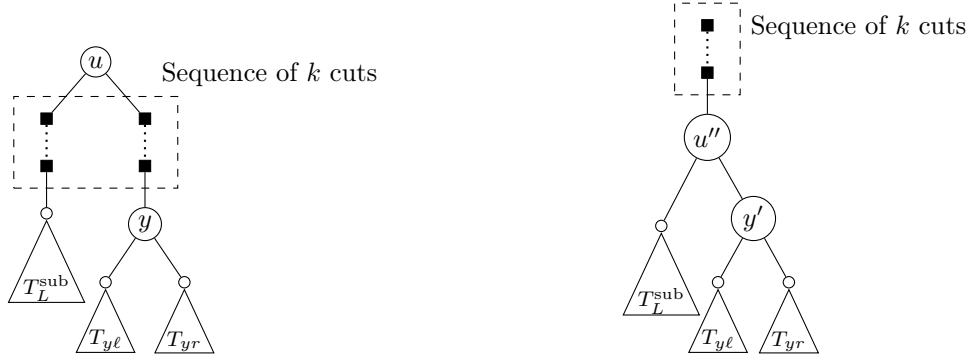
\ifspringer
\else
\par
\fi
    Observe that 
        \begin{equation*}
        \begin{aligned}
            \treetime_{\W}(T_{u}) 
            &= \W(u) + 2 \sum_{i \in [k]} \W(v_i) + \W(y) + \treetime_{\W}(T_{y\lbd}) + \treetime_{\W}(T_{y\rbd}) + \treetime_{\W}(T^{\mathrm{sub}}_L),
            \\
        \treetime_{\W}(T_{u}^{\prime\prime})
        &= \sum_{i \in [k]} \W(v_{i}^{\prime\prime}) + \W(u^{\prime\prime}) + \W(y^\prime) + \treetime_{\W}(T_{y\lbd}) + \treetime_{\W}(T_{y\rbd}) + \treetime_{\W}(T^{\mathrm{sub}}_L)\\
        &= \sum_{i \in [k]} \W(v_{i}) + 2\W(y) + \treetime_{\W}(T_{y\lbd}) + \treetime_{\W}(T_{y\rbd}) + \treetime_{\W}(T^{\mathrm{sub}}_L).
        \end{aligned}
    \end{equation*}
    Then 
        $
            \treetime_{\W}(T_{u}^{\prime\prime})-\treetime_{\W}(T_{u}) 
            = \W(y) - 
            ( \W(u) + \sum_{i \in [k]} \W(v_i) )
        $
    is always nonpositive by \eqref{equa:root2}.
    Hence, the modified tree $T^{\prime\prime}$ has tree time at most that of $T$,
    while having $k$ fewer cut nodes,
    contradicting the assumption that $T$ is cut-minimal.
\end{proof}

We next prove, in \Cref{thm:rootcuts}, that under symmetric branching, every cut-minimal tree has all cuts at the root, allowing multiple cut families and cut generation time.
This generalizes a previous result~\cite[Theorem~23]{KazLeBSan24_abstractbc-mpb}
that assumes only one cut family and constant cut strength.

\begin{theorem}\label{thm:rootcuts}
    Let $(\vec{\cbd},\lbd,\rbd) \in \cutfunset^\ncutfams \times \posreals \times \posreals$
    such that $\lbd = \rbd$.
    For any target bound $\targetbound$ and generalized time function $\W$,
    every cut-minimal \bc tree proving bound $\targetbound$ has only root cuts.
\end{theorem}
\begin{proof}
    Let $T$ be a cut-minimal tree proving bound $\targetbound$ with respect to $\W$, and suppose that a branch node $u$ of $T$ has a cut node $v$ as a child. Let $y$ be the other child of $u$.
    Since $\lbd = \rbd$, $z_T(v) = z_T(y)$ and $\historyvec(v) = \historyvec(y)$, so (as $T$ is cut-minimal) $T_v$ and $T_y$ have the same tree time and contain the same number of cut nodes.
    Replacing $T_y$ by a copy of $T_v$ therefore yields a cut-minimal tree in which both child subtrees of $u$ begin with the same maximal sequence of $k \ge 1$ cut nodes, with the same families in the same order, followed by a branch or leaf node, contradicting \Cref{lem:nocommonprefix}.
    Hence, no branch node of $T$ has a cut node as a child, so every cut node of $T$ lies on the path of cut nodes starting at the root, i.e., $T$ has only root cuts.
\end{proof}

We emphasize that symmetry plays a crucial role in this proof. In particular, symmetry of a tree will force cut nodes to move up, regardless of how the time function is defined. 
On the other hand, the proof does not rely on branching strength: the result can be extended even with multiple symmetric branch variable candidates, i.e., possibly different values throughout the tree but always equal left/right bound improvements.

\subsection{Ordering Weak Versus Strong Cut Families}
\label{sec:sparse-then-dense}

Restricting to only root-cut trees, whether optimal or not, leaves open questions such as the number of cuts to add (previously approximated in \cite[Theorem~15]{KazLeBSan24_abstractbc-mpb}) and, with multiple cut families, their order.
We consider this latter question for a path of cuts (which may not be at the root) in a minimal-time tree, when there are weak and strong (constant-strength) cut families. %
For a constant cut function, we use $\cbd_i$ to denote both the function and its value.

\begin{restatable}{theorem}{sparsethendense}
\label{thm:sparse-then-dense}
    Let $(\vec{\cbd},\lbd,\rbd) \in \cutfunset^\ncutfams \times \posreals \times \posreals$ with constant cut functions $\cbd_i$ for $i \in [\ncutfams]$.
    For any target bound $\targetbd$ and generalized time function $\W$, there exists a minimal-time tree proving $\targetbd$ in which, along every path of cut nodes, the cut families appear in nondecreasing order of strength $\cbd_i$.
\end{restatable}

\begin{proof}
    Let $T$ be a minimal-time tree proving $\targetbd$ with respect to $\W$, and consider a path of $\kappa$ consecutive cut nodes in $T$.
    Let $k_i$ denote the number of cuts from family $i \in [\ncutfams]$ in this path.
    We will reorder these cut nodes by nondecreasing strength.
    Since the cut functions are constant, the reordered tree still proves $\targetbd$ (any permutation of the path increases the bound by $\sum_{i=1}^{\ncutfams} k_i \cbd_i$), and all nodes outside of the path keep their bounds and node times.
    For the nodes on the path, as $\W$ is the sum of a simple time function $\w$ and separation times, we analyze its two components separately.

    First, we show that total separation time is invariant under any reordering.
    Let $\historyvec$ be the history vector at the first cut node in the path.
    The total separation time of the sequence is 
        $\sum_{i=1}^{\ncutfams} \sum_{h=\historyvec_i}^{\historyvec_i+k_i-1} \septime_i(h)$, 
    which depends only on $\historyvec$ and the counts $k_i$, not the order.

        Next, for node processing time, label the positions of the cuts as $0, \dots, \kappa-1$, from start to end of the path, and let $f_j$ be the index of the cut family at position $j$. Let $v_{j_0}$, at position $j_0$, be a cut from a weakest family in the sequence, and consider moving $v_{j_0}$ to position $0$. This does not affect the nodes before position $0$ or after position $j_0$, so it suffices to analyze positions $0, \dots, j_0$. Let $\w_j$ and $\w'_j$ denote the node processing time at position $j$ before and after the move. Then $\w_0 = \w'_0$. 
        For $j \in [j_0]$, let $X_j \defeq \historyvec + \sum_{j'=0}^{j-2} e_{f_{j'}}$, from which
        the history vector at position $j$ is $X_j + e_{f_{j-1}}$ before the move and $X_j + e_{f_{j_0}}$ after it. Since $\cbd_{f_{j_0}} \le \cbd_{f_{j-1}}$, \Cref{def:simple-time-function} gives $\w'_j = \w(X_j + e_{f_{j_0}}) \le \w(X_j + e_{f_{j-1}}) = \w_j$.

        Recursively repeating this move of the next-weakest cut node, starting at position $1$, then $2$, and so on, orders the path by nondecreasing strength without increasing its tree time.
        Applying this to every maximal path of consecutive cut nodes in $T$ (these sequences are disjoint, and each step leaves all other nodes unchanged) yields a minimal-time tree with the desired property.
\end{proof}

\begin{figure}[t]
    \centering
    \begin{minipage}[t]{0.49\textwidth}
        \centering
        \begin{forest}
        for tree = {small node,l=1pt, s sep= 25pt}
            [,cut, fill = none, lvlwt = {${\w(0,0)}+\septime_1(0)$},  label={[fill=none,font=\small]180:{{weak}}},
                        [,cut,lvlwt = {${\w(1,0)}+\septime_2(0)$},  label={[fill=none,font=\small]180:{{strong}}},
                            [,cut, lvlwt = {${\w(1,1)}+\septime_2(1)+M$},  label={[fill=none,font=\small]180:{{strong}}}, 
                            [,lvlwt = {${\w(1,2)}$},fill = none, inner sep = 0.5pt, 
                            ]
                            ]
                        ]
            ]
    \end{forest}
    \end{minipage}
    \hspace{1.5em}
    \begin{minipage}[t]{0.45\textwidth}
        \centering
        \begin{forest}
        for tree = {small node,l=1pt, s sep= 25pt}
            [,cut, lvlwt = {${\w(0,0)}+\septime_2(0)$},
                        [,cut,fill = none, lvlwt = {${\w(0,1)}+\septime_1(0)$},
                            [,cut, lvlwt = {${\w(1,1)}+\septime_2(1)$},
                            [,lvlwt = {${\w(1,2)}$},fill = none, inner sep = 0.5pt, 
                            ]
                            ]
                        ]
            ]
    \end{forest}
    \end{minipage}
    \caption{
        With more general separation time functions,
        placing weaker cuts before stronger ones, as in \Cref{thm:sparse-then-dense}, may be suboptimal.
        In the example in this figure, if generating two consecutive cuts from the same family incurs a cost $M > \w(0,1) - \w(1,0)$,
        the alternating structure in the right panel is cheaper than the weak-then-strong pattern in the left panel.
        }
    \label{fig:sparse-then-dense}
\end{figure}
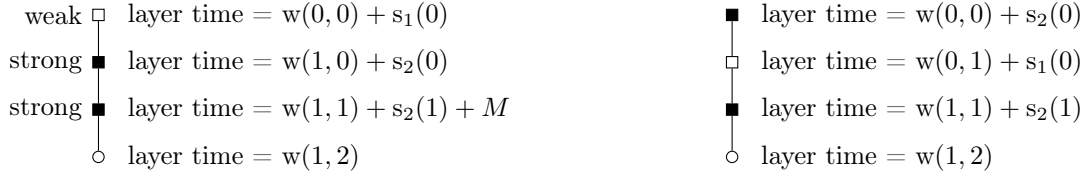

\Cref{thm:sparse-then-dense} relies on the crucial assumption from \Cref{def:separator-time-function} that separation time depends only on the number of previous cuts from the same family, so that separation time does not affect the order of consecutive cut nodes.
\Cref{fig:sparse-then-dense} depicts one example where this assumption is relaxed:
when there is a heavy penalty on consecutive nodes belonging to the same cut family,
the minimal-time tree alternates between cut families instead of following a weak-then-strong order.

\subsection{Root-Cut Trees Under Asymmetric Branching}
\label{sec:124}

In \Cref{sec:root}, we showed that it is sufficient to consider root cuts in symmetric branching settings with multiple cut families.
We now turn to the single cut-family, single branching-family \bc model with constant but nonidentical bound improvements ($\lbd \neq \rbd$), under a simple time function.

By \Cref{lem:nocommonprefix}, symmetry in a tree shifts cut nodes upward toward the root.
Consequently, a natural question is whether symmetry in the tree is necessary, i.e., might one always find, under an asymmetric branching regime, a combination of target bound and time function that forces the inclusion of local cuts within any minimal-time tree.
We prove that this is not universally true: there exists an asymmetric setting 
    ($\cbd\mathbin{:}\lbd\mathbin{:}\rbd$ in ratio $1\mathbin{:}2\mathbin{:}4$) 
in which a minimal-time tree without any internal cuts always exists, regardless of the specific form of the (simple) time function.
Proving this involves several intermediate results that we present first. %

Define a \emph{branch-then-bothcut} node to be a branch node with both children being cut nodes from the same family.
The following corollary, which follows directly from \cref{lem:nocommonprefix}, states that a branch-then-bothcut node cannot exist for a single cut family setting.%
\footnote{It can be optimal for both children of a branch node to be cuts from different families.}

\begin{restatable}{corollary}{nobothcut}
\label{cor:nobothcut}
     For any $(\cbd,\lbd,\rbd) \in \posreals^3$, target bound $\targetbd$, and generalized time function $\W$,
     every cut-minimal tree proving $\targetbd$ has no branch-then-bothcut node.
\end{restatable}
\begin{proof}
    If a cut-minimal tree has a branch-then-bothcut node, then the two children of this node
    begin with the same sequence of $\min\{k_\lbd, k_\rbd\}$ cut nodes, where $k_\lbd$ and $k_\rbd$ are the lengths of the maximal cut paths starting at each child.
    By maximality of the cut paths, at least one of these sequences is followed by a branch or leaf node,
    so \Cref{lem:nocommonprefix} implies the tree is not cut-minimal.
\end{proof}

It remains unclear whether branch nodes with exactly one child as a cut node may be necessary. More importantly, if such branch nodes do appear in a minimal-time tree, it is not immediately clear whether the deepest such node has a left-child cut or a right-child cut.
We call a branch-then-cut node a \emph{branch-then-leftcut} node if its left child is an internal cut.
We address this question (for internal cut nodes) in \Cref{lem:124deepest}.
When $\lbd = \rbd$, \Cref{thm:rootcuts} implies that no cut-minimal tree has a branch-then-cut node, so we focus on $\rbd = 2\lbd$.

\begin{restatable}{lemma}{onetwofourdeepest}
\label{lem:124deepest}
    Let $(\cbd,\lbd,\rbd) \in \posreals^3$ with $\rbd = 2\lbd$,
    $\w$ be a simple time function, and $\targetbd$ be a target bound.
    Among all cut-minimal trees proving $\targetbd$, let $T$ be one whose deepest branch-then-cut node (if any) is as deep as possible.
    Every deepest branch-then-cut node of $T$ is a branch-then-leftcut node.
\end{restatable}

\begin{proof}
    Suppose that a deepest branch-then-cut node $u$ of $T$ has an internal cut node $v$ as its right child.
    The left child of $u$ is a branch node $x$
    --- it is not a cut node, as otherwise $u$ would be a branch-then-bothcut node, contradicting \Cref{cor:nobothcut}, and it is not a leaf node, as $z_{T}(x) < z_{T}(v) < \targetbd$,
    where the last inequality holds by minimality since $v$ is not a leaf.
    Let $y$ be the left child of $x$.
    Then $z_{T}(y) = z_T(u) + 2\lbd = z_{T}(v)$, and $v$ and $y$ have the same cut history, so either of the subtrees $T_v$ and $T_y$ can replace the other, and minimality gives $\treetime_\w(T_v) = \treetime_\w(T_y)$.
    As $T$ is cut-minimal, $T_y$ has at least as many cut nodes as $T_v$, so replacing $T_y$ by a copy of $T_v$ yields a cut-minimal tree in which $x$ is a branch-then-cut node deeper than $u$, contradicting the choice of $T$.
\end{proof}

Unlike symmetric branching, where symmetry is maintained throughout the tree, any node-shifting operation can significantly alter the tree structure under an asymmetric branching setting. Therefore, it is crucial to understand how tree structure changes as the target bound varies. 
\Cref{lem:124structure} provides a useful pattern under the special setting of $\cbd \mathbin{:} \lbd \mathbin{:} \rbd = 1 \mathbin{:} 2 \mathbin{:} 4$.
For any $z \in \R$, let $p(z) \defeq \cardinality{\PureBB{z}}$ (the size of the smallest pure branching tree that proves bound $z$), with $p(z) = 1$ if $z \le 0$, and let $\nleaf(z)$ denote the number of leaf nodes of $\PureBB{z}$ with the smallest bound.
If $\lbd$, $\rbd$, and $z$ are even,
then $\nleaf(z)$ equals the number of leaf nodes proving exactly $z$. 
Note also that $p(0) = \nleaf(0) = 1$.

\begin{lemma}\label{lem:124structure}
    If $\cbd \mathbin{:} \lbd \mathbin{:} \rbd = 1 \mathbin{:} 2 \mathbin{:} 4$ and $\targetbd \ge 4\cbd$ is an even multiple of $\cbd$,
    then
    \begin{enumerate}[(1)]
        \item $p(\targetbd + 2\cbd) + p(\targetbd - 4\cbd) = 2p(\targetbd)$, and
        \item $\nleaf(\targetbd + 2\cbd) + \nleaf(\targetbd - 4\cbd) = 2\nleaf(\targetbd)$, with $\nleaf(\targetbd) - \nleaf(\targetbd - 4\cbd) = \nleaf(\targetbd - 2\cbd) \ge 1$.
    \end{enumerate}
\end{lemma}
\begin{proof}
    For $z > 0$, the root of $\PureBB{z}$ branches into copies of $\PureBB{z - 2\cbd}$ and $\PureBB{z - 4\cbd}$, so $p(z) = 1 + p(z - 2\cbd) + p(z - 4\cbd)$~\cite{LeBNem17}. Moreover, for $z \ge 4\cbd$, the leaves of $\PureBB{z}$ that prove exactly $z$ are the union of the leaves of $\PureBB{z - 2\cbd}$ and $\PureBB{z - 4\cbd}$ that prove their respective (relative) target bounds,
    so $\nleaf(z) = \nleaf(z - 2\cbd) + \nleaf(z - 4\cbd)$.
    Hence, $f \in \{p, \nleaf\}$ satisfies $f(z) = a_f + f(z - 2\cbd) + f(z - 4\cbd)$ for $z \ge 4\cbd$, with $a_p = 1$ and $a_{\nleaf} = 0$.
    This gives $f(\targetbd + 2\cbd) - f(\targetbd) = f(\targetbd) - f(\targetbd - 4\cbd)$, which proves (1) and the first part of (2).
    Finally, 
        $\nleaf(\targetbd) - \nleaf(\targetbd - 4\cbd)
        = \nleaf(\targetbd - 2\cbd) \ge 1$,
    as the leftmost leaf of $\PureBB{\targetbd - 2\cbd}$ proves exactly $\targetbd - 2\cbd$.
\end{proof}

Next, \Cref{lem:oneleafcut} states that, when $\lbd \ge 2\cbd$, it is (weakly) better to replace two consecutive leaf cuts (two cut nodes followed by a leaf) by a branch node with two leaf children, i.e., no leaf is preceded by two consecutive cut nodes.

\begin{lemma}\label{lem:oneleafcut}
    For any $(\cbd,\lbd,\rbd) \in \posreals^3$ with 
        $2 \cbd \le \lbd \le \rbd$,
    simple time function $\w$, and target bound $\targetbd$, every cut-minimal tree proving $\targetbd$ contains at most one leaf cut node on each root-to-leaf path.
\end{lemma}

\begin{proof}
    Suppose that a cut-minimal tree contains two consecutive cut nodes followed by a leaf, with $h$ cut nodes above the first of them.
    Replacing these three nodes by a branch node with two leaf children yields a tree that 
        still proves $\targetbd$ (as $\rbd \ge \lbd \ge 2\cbd$), 
        has weakly smaller tree time (as $3\w(h) \le \w(h) + \w(h+1) + \w(h+2)$), 
    and has two fewer cut nodes, a contradiction.
\end{proof}

\Cref{lem:replace-leaf-cuts-by-pure-branching} proves a property of a cut-minimal tree
when $\cbd \mathbin{:} \lbd \mathbin{:} \rbd = 1\mathbin{:}2\mathbin{:}4$:
for a subtree $T_g$ with no internal cuts,
after we ``remove'' any leaf cuts (i.e., replace each leaf cut and its leaf child by a single leaf),
we obtain the smallest pure branching tree $T_G$ proving a bound $\targetbd_G$ that is an even multiple of $\cbd$.
Further, either $T_g = T_G$, or $T_g$ is exactly $T_G$ with an additional cut node inserted above 
each leaf $v$ of $T_G$ with $z_{T_g}(v) = \targetbd_G$.
Thus, the time of $T_g$ is determined by $\targetbd_G$ and whether leaf cuts are used.
The proof of \Cref{thm:124} will use an idea similar to our proof in \cref{lem:nocommonprefix}, of shifting cuts up or down a candidate tree.
Our analysis is complicated by the asymmetric setting; \Cref{lem:124structure} handles this, as our moves change $\targetbd_G$ by even multiples of $\cbd$.

\begin{lemma}
\label{lem:replace-leaf-cuts-by-pure-branching}
    Suppose that $\cbd \mathbin{:} \lbd \mathbin{:} \rbd = 1\mathbin{:}2\mathbin{:}4$, and let $T$ be a cut-minimal tree proving a target bound $\targetbd$ with respect to a simple time function $\w$.
    Let $g$ be a branch node of $T$ such that the subtree $T_g$ contains no internal cuts, 
    let $T_G$ be the pure branching tree obtained from $T_g$ by removing its leaf cuts, and 
    let $\targetbd_G$ be the smallest bound at a leaf node of $T_G$, relative to $g$.
    Then
    \begin{enumerate}[(1)]
        \item $\targetbd_G$ is an even multiple of $\cbd$ with $\targetbd_G \ge 2\cbd$;
        \item $T_G = \PureBB{\targetbd_G}$, the smallest pure branching tree proving $\targetbd_G$; and
        \item the nodes $v$ of $T_g$ with $z_{T_g}(v) = \targetbd_G$ are either all leaves or all leaf cuts, and $z_{T_g}(v) = \targetbd_G$ at every leaf cut $v$ in $T_g$.
    \end{enumerate}
\end{lemma}

\begin{proof}
    As $(\lbd,\rbd) = (2\cbd,4\cbd)$, for all nodes $v$ of $T_G$, $z_{T_G}(v)$ is an even multiple of $\cbd$. Since $g$ is a branch node, we have $\targetbd_G \ge 2\cbd$, proving (1).
        Let $R \defeq \targetbd - z_T(g)$ be the remaining bound at $g$,
        so that $z_{T_g}(v) \ge R$ for every leaf $v$ of $T_g$.
    Let $v$ be a leaf node of $T_G$.
    By minimality, if $z_{T_G}(v) \ge R$, then $T_g$ has the exact same leaf (without a leaf cut inserted).
    Conversely, if $z_{T_G}(v) < R$, then a leaf cut is needed there in $T_g$,
    and, by \Cref{lem:oneleafcut}, a single leaf cut suffices.
    Thus, as $\targetbd_G$ is the smallest leaf bound in $T_G$, $T_g$ contains leaf cuts if and only if $\targetbd_G < R$.
    Moreover, $R \le \targetbd_G + \cbd$, and, as all bounds are even multiples of $\cbd$, every leaf cut $u$ of $T_g$ has $z_{T_g}(u) = \targetbd_G$, proving the second half of (3).

    \looseness=-1
    Observe that, by construction, all branch nodes in $T_G$ and $T_g$ are identical.
    If $T_g$ contains no leaf cuts,
    then $T_g = T_G$ and $R \le \targetbd_G$, so (2) follows from minimality of $T$, and all nodes of $T_g$ with bound $\targetbd_G$ are leaves, completing (3).
    Otherwise, $\targetbd_G < R$, and
    let $v$ be a branch node of $T_G$ with bound $z_v \defeq z_{T_G}(v)$ ($= z_{T_g}(v)$).
    For $T_G$ to be the smallest pure branching tree proving $\targetbd_G$, we need to show that $z_v < \targetbd_G$.
    By minimality of $T$, as $v$ is not a leaf, $z_v < R$ (implying $z_v \le \targetbd_G$).
    By the definition of $\targetbd_G$, some leaf of $T_G$ has bound $\targetbd_G < R$, and, as argued above,
    the corresponding position in $T_g$ is a leaf cut.
    If $z_v = \targetbd_G$, the children of $v$ are leaves with bound at least $R$, and $T_v$ can replace the leaf cut and leaf in $T_g$ without increasing tree time, contradicting cut-minimality of $T$.
    Thus, $z_v < \targetbd_G$, proving (2).
    Finally, as $\targetbd_G < R$, every leaf $u$ of $T_G$ with $z_{T_G}(u) = \targetbd_G$ needs a leaf cut in $T_g$, and, by (2), these are the only nodes of $T_g$ with bound $\targetbd_G$, proving the first half of (3).
\end{proof}

Finally, we prove  \Cref{thm:124}, the main result of this subsection:
for $(\cbd,\lbd,\rbd)$ satisfying a ratio $1\mathbin{:}2\mathbin{:}4$, there exists a minimal-time tree where all cuts are either at the root or at the leaves.

\begin{restatable}{theorem}{onetwofour}
\label{thm:124}
    For $(\cbd,\lbd,\rbd) \in \posreals^3$ with $\cbd\mathbin{:}\lbd\mathbin{:}\rbd = 1\mathbin{:}2\mathbin{:}4$, every target bound $\targetbd$, and any simple time function $\w$, 
    there exists a minimal-time tree that does not contain any internal cuts.
\end{restatable}
\begin{proof}
    \looseness=-1
    Suppose, for the sake of contradiction, that for some target bound $\targetbd$ and time function~$\w$, every minimal-time tree proving $\targetbd$ has an internal cut.
    Then every cut-minimal tree has a branch-then-cut node.
    Let $u$ be the deepest branch-then-cut node of a tree $\hat{T}$ whose deepest branch-then-cut node is as deep as possible among all cut-minimal trees proving $\targetbd$.
    Let $T$ be the subtree rooted at $u$, with $T_\ell$ the left and $T_\rbd$ the right child subtree of $u$. 
    By \Cref{lem:124deepest}, $u$ is a branch-then-leftcut node.
    Let $\kappa$ be the length of the maximal path of cut nodes starting at the root of $T_\ell$.
    These cuts are internal, so this path is followed by a branch node, which we will denote by $g$, the root of subtree $T_g$.

    \looseness=-1
    The proof compares two ways of moving this path: up a layer above $u$, or down two layers below two branch nodes.
    As pure branching trees have a recursive structure in this setting (\Cref{lem:124structure}), these moves change tree time by exactly opposite amounts.
    As moving the cuts down cannot help in a minimal-time tree, moving them up is at least as good. This yields a tree with fewer cut nodes or a deeper branch-then-cut node, contradicting the choice of $\hat{T}$.

    Let $K$ be the number of cut nodes above $u$, and let $\bar\w(h) \defeq \w(K+h)/\w(K)$ for $h \in \nonnegints$, which is again a simple time function with $\bar\w(0) = 1$.
    As the roots of $T$, $T_\ell$, and $T_\rbd$ all have $K$ cut nodes above them, \cite[Lemma~24]{KazLeBSan24_abstractbc-mpb} implies that they are minimal-time trees for their remaining bounds with respect to $\bar\w$; for simplicity, we write $\w$ for $\bar\w$ hereafter.
    Moreover, these trees are also cut-minimal for their remaining bounds, as otherwise replacing one of them by a minimal-time subtree with fewer cut nodes would contradict the cut-minimality of the whole tree.

    Let $\targetbd_u \defeq \targetbd - z_{\hat T}(u)$ be the remaining bound at $u$, and let $\targetbd_\ell \defeq \targetbd_u - \lbd$ be the remaining bound that $T_\ell$ needs to prove.
    By minimality and the need for an internal cut in $T_\ell$, we have that $\targetbd_\ell > \lbd = 2\cbd$.
    This implies that $T_\rbd$ is not a single leaf node,
    as it must prove $\targetbd_\ell - 2\cbd > 0$.
    The root of $T_\rbd$ is not a cut node, as otherwise $u$ would be a branch-then-bothcut node, contradicting \Cref{cor:nobothcut}.
    Thus, the root of $T_\rbd$ is a branch node. 
    Let $T_1$ and $T_2$ be its left and right child subtrees.

    Since $u$ is a deepest branch-then-cut node, there is no internal cut in $T_\rbd$ or in $T_g$.
    Let $T_G$ and $\targetbd_G$ be as in \Cref{lem:replace-leaf-cuts-by-pure-branching} applied to $g$, and let $I_{\text{leaf}} \in \{0,1\}$ indicate whether $T_g$ contains a leaf cut.
    By \Cref{lem:replace-leaf-cuts-by-pure-branching}, $\targetbd_G \ge 2\cbd$ is an even multiple of $\cbd$ and $T_G = \PureBB{\targetbd_G}$ (statements (1) and (2)), and $T_g$ is obtained from $T_G$ by inserting a leaf cut above each of its $\nleaf(\targetbd_G)$ leaves $v$ with $z_{T_g}(v) = \targetbd_G$ if $I_{\text{leaf}} = 1$ (statement (3)).
    Thus, $T_\ell$ proves $\kappa\cbd + \targetbd_G + I_{\text{leaf}}\cbd \ge \targetbd_\ell$.

    For an even multiple $\bar\targetbd \ge 0$ of $\cbd$, define $F(\bar\targetbd,I_{\text{leaf}}) \defeq \w(\kappa)p(\bar\targetbd) + I_{\mathrm{leaf}}\w(\kappa+1)\nleaf(\bar\targetbd)$,
    so that $F(\targetbd_G, I_{\text{leaf}}) = \treetime_{\w}(T_g)$.
    The tree times of $T_\ell$, $T_\rbd$, and $T$ are then
    
  \begin{equation}\label{eq:124-Tl}
   \begin{aligned}
       \tau_\w(T_\ell) &= \sum_{i\in [0,\kappa-1]} \w(i) + F(\targetbd_G, I_{\text{leaf}}), 
       \\ 
       \tau_\w(T_\rbd) &= 1 + \tau_\w(T_1) + \tau_\w(T_2), 
       \quad\mbox{and}\quad
       \tau_\w(T) = 1 + \tau_\w(T_\ell) + \tau_\w(T_\rbd).
   \end{aligned}
   \end{equation}
    As $F$ is a linear combination of $p$ and $\nleaf$, \Cref{lem:124structure} implies that the identity
    \begin{equation}\label{eq:124-F}
        F(\bar\targetbd+2\cbd, I_{\text{leaf}}) - F(\bar\targetbd,I_{\text{leaf}}) = F(\bar\targetbd,I_{\text{leaf}}) - F(\bar\targetbd-4\cbd,I_{\text{leaf}})
    \end{equation}
    holds for all even multiples $\bar\targetbd \ge 4\cbd$ of $\cbd$ and $I_{\text{leaf}} \in \{0,1\}$.

\begin{figure}
  \centering
  \begin{minipage}[t]{0.3\textwidth}\centering
    \begin{forest} baseline, for tree={small node, l=8pt, s sep=14pt}
    [$u$, fill=none, inner sep=1.5pt, s sep=34pt
      [, cut, name=TLroot [, cut, cutpath [, fill=none, tri={T_g}{11mm}, name=TLg]]]
      [, fill=none, name=TRroot, s sep=22pt
        [, fill=none, tri={T_1}{7mm}, name=TRl]
        [, fill=none, tri={T_2}{7mm}, name=TRr]]]
    \draw[dashed, rounded corners] ([xshift=-4.6mm,yshift=-8.4mm]TRl.south) rectangle ([xshift=4.6mm,yshift=8.8mm]TRr.south);
    \node[anchor=south east, font=\small, inner sep=1.5pt] at ([xshift=4.6mm,yshift=8.8mm]TRr.south) {$T_\rbd$};
    \draw[dashed, rounded corners] ([xshift=-6.2mm,yshift=-12.4mm]TLg.south) rectangle ([xshift=6.2mm,yshift=1mm]TLroot.north);
    \node[anchor=south west, font=\small, inner sep=1.5pt] at ([xshift=-6.2mm,yshift=1mm]TLroot.north) {$T_\lbd$};
    \end{forest}
  \end{minipage}\hfill
  \begin{minipage}[t]{0.38\textwidth}\centering
    \begin{forest} baseline, for tree={small node, l=8pt, s sep=16pt}
    [$u$, fill=none, inner sep=1.5pt, s sep=50pt
      [, fill=none, s sep=26pt
        [, fill=none, s sep=22pt
          [, cut [, cut, cutpath [, fill=none, tri={T'_g}{7mm}, name=TLpg]]]
          [, fill=none, tri={T_2}{7mm}]]
        [, fill=none, tri={T_1}{7mm}, name=TLp1]]
      [, fill=none, tri={T_\rbd}{10mm}]]
    \coordinate (TLpBL) at ([xshift=-4.6mm,yshift=-8.4mm]TLpg.south);
    \coordinate (TLpTR) at ([xshift=4.6mm,yshift=8.8mm]TLp1.south);
    \draw[dashed, rounded corners] (TLpBL) rectangle (TLpTR);
    \node[anchor=south west, font=\small, inner sep=1.5pt] at (TLpBL |- TLpTR) {$T'_\lbd$};
    \end{forest}
  \end{minipage}\hfill
  \begin{minipage}[t]{0.22\textwidth}\centering
    \begin{forest} baseline, for tree={small node, l=8pt}
    [, cut, label={[font=\small]180:{$u$}}
      [, cut, cutpath [, fill=none, tri={T''_g}{14mm}]]]
    \end{forest}
  \end{minipage}
  \caption{Left panel: tree $T$.
    Middle panel: tree $T'$, obtained by moving the $\kappa$ cut nodes down two layers.
    Right panel: tree $T''$, obtained by moving the $\kappa$ cut nodes up one layer to $u$.}
  \label{fig:124}
\end{figure}
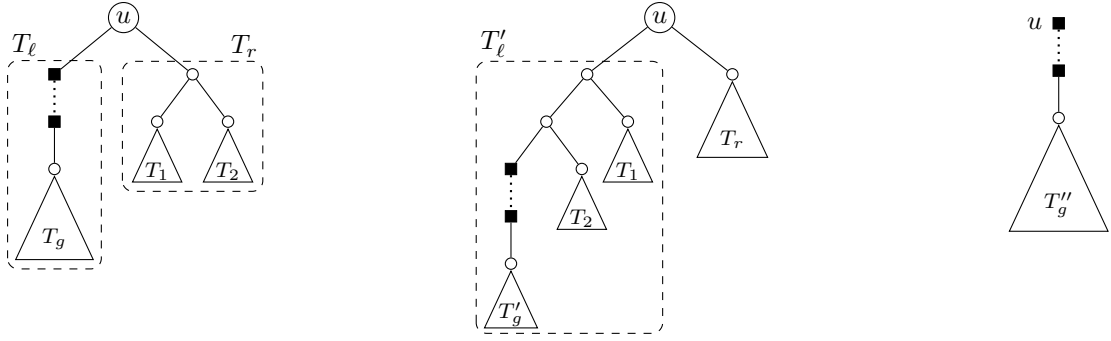

    \mbox{}\\[-1ex]\noindent
        \emph{Case 1: $\targetbd_G \ge 4\cbd$.}        
    We compare $T$ with two other trees, illustrated in \Cref{fig:124}, that also prove $\targetbd_u$.
    For the first tree, $T'$, we shift the $\kappa$ cut nodes \emph{down two layers} along the leftmost path of $T_{\ell}$ while keeping $T_{r}$ unchanged.
    In $T'$, the left child of $u$ is a branch node whose right child is a copy of $T_1$ (as both are reached from $u$ by one left and one right branch, they have the same bound and history).
    The left child of this branch node is again a branch node, whose right child is a copy of $T_2$ (as branching twice left and once right is equivalent to two right branches when $\rbd = 2 \lbd$),
    and whose left child is the path of $\kappa$ cut nodes followed by $T'_g$,
    obtained from
    $T'_G \defeq \PureBB{\targetbd_G - 4\cbd}$ by inserting (if $I_{\text{leaf}} = 1$) a leaf cut above each leaf $v$ with $z_{T'_G}(v) = \targetbd_G - 4\cbd$.
    To see that $T'$ proves $\targetbd_u$, the only nontrivial case is the leftmost path to the leaves of $T'_g$.
    As all bounds are even multiples of $\cbd$, every leaf $v$ of $T'_G$ with $z_{T'_G}(v) \neq \targetbd_G - 4\cbd$ has $z_{T'_G}(v) \ge \targetbd_G - 2\cbd$, so $T'_g$ proves $\targetbd_G - 4\cbd + I_{\text{leaf}}\cbd$.
    Hence, every leaf $v$ of $T'$ in $T'_g$ satisfies 
        $z_{T'}(v) 
        \ge 3\lbd + \kappa\cbd + (\targetbd_G - 4\cbd) + I_{\text{leaf}}\cbd 
        = \lbd + \kappa\cbd + \targetbd_G + I_{\text{leaf}}\cbd 
        \ge \lbd + \targetbd_\ell
        = \targetbd_u$.
    
                Next, define tree $T''$ by shifting the $\kappa$ cut nodes \emph{up one layer}: 
                $T''$ has a path of $\kappa$ cut nodes followed by $T''_g$, analogously related to $T''_G \defeq \PureBB{\targetbd_G + 2\cbd}$ as $T'_g$ to $T'_G$.
                    As with $T'_g$, every leaf $v$ of $T''$ satisfies $z_{T''}(v) \ge \kappa\cbd + (\targetbd_G + 2\cbd) + I_{\text{leaf}}\cbd = \lbd + \kappa\cbd + \targetbd_G + I_{\text{leaf}}\cbd \ge \targetbd_u$, so $T''$ also proves $\targetbd_u$.

    The subtree $T'_\ell$ of $T'$ rooted at the left child of $u$ consists of two branch nodes, copies of $T_1$ and $T_2$, the $\kappa$ cut nodes, and $T'_g$, so $\tau_\w(T'_\ell) = 1 + \tau_\w(T_\rbd) + \sum_{i=0}^{\kappa-1}\w(i) + F(\targetbd_G-4\cbd, I_{\text{leaf}})$, using $\tau_\w(T_\rbd) = 1 + \tau_\w(T_1) + \tau_\w(T_2)$.
                The tree time of $T''$ is $\tau_\w(T'') = \sum_{i=0}^{\kappa-1}\w(i) + F(\targetbd_G+2\cbd, I_{\text{leaf}})$.
    Combining these expressions with \Cref{eq:124-Tl,eq:124-F} gives us
    \begin{equation}\label{eq:124-key}
    \begin{aligned}
        \Delta \coloneqq \tau_\w(T) - \tau_\w(T'')
        &= 1 + \tau_\w(T_\rbd) - \big(F(\targetbd_G + 2\cbd, I_{\text{leaf}}) - F(\targetbd_G, I_{\text{leaf}})\big)\\
        &= 1 + \tau_\w(T_\rbd) - \big(F(\targetbd_G, I_{\text{leaf}}) - F(\targetbd_G - 4\cbd, I_{\text{leaf}})\big) \\
        &= \tau_\w(T'_\ell) - \tau_\w(T_\ell) = \tau_\w(T') - \tau_\w(T).
    \end{aligned}
    \end{equation}

    Since $T$ is minimal-time and both $T'$ and $T''$ prove $\targetbd_u$, we have $\tau_\w(T') \ge \tau_\w(T)$ and $\tau_\w(T'') \ge \tau_\w(T)$, i.e., $\Delta \ge 0$ and $\Delta \le 0$ by \Cref{eq:124-key}.
    Thus, $\Delta = 0$, and $T'$ and $T''$ are both minimal-time trees proving $\targetbd_u$.
    Let $n(\cdot)$ denote the number of cut nodes of a tree.
    As $T'$ contains $T_\rbd$ as well as copies of $T_1$ and $T_2$, statement~(2) of \Cref{lem:124structure} gives
    \begin{equation*}
    \begin{aligned}
        n(T') + n(T'')
        &= \big(\kappa + I_{\text{leaf}}\,\nleaf(\targetbd_G - 4\cbd) + 2n(T_\rbd)\big) + \big(\kappa + I_{\text{leaf}}\,\nleaf(\targetbd_G + 2\cbd)\big) \\
        &= 2\big(\kappa + I_{\text{leaf}}\,\nleaf(\targetbd_G) + n(T_\rbd)\big)
        = 2n(T).
    \end{aligned}
    \end{equation*}
    As $T$ is cut-minimal, $n(T') \ge n(T)$ and $n(T'') \ge n(T)$, so $n(T') = n(T'') = n(T)$, and replacing $T$ by $T'$ in $\hat T$ yields a cut-minimal tree.
    If $\targetbd_G \ge 6\cbd$, then the $\kappa$ cut nodes of $T'$ are followed by a branch node, so this cut-minimal tree has a branch-then-cut node two layers deeper than $u$, contradicting the choice of $\hat T$.
    If $\targetbd_G = 4\cbd$, then $T'_G = \PureBB{0}$ is a single leaf, so in $T'$ the $\kappa$ cut nodes, together with a leaf cut if $I_{\text{leaf}} = 1$, are followed by a leaf, which contradicts \Cref{lem:oneleafcut} unless $\kappa = 1$ and $I_{\text{leaf}} = 0$.
    
    In this remaining case, $1 = n(T'') = n(T) = 1 + n(T_\rbd)$, so $T_\rbd$ has no cut nodes.
    Then $T_\ell$ consists of a cut node followed by $T_G = \PureBB{4\cbd}$, which has $5$ nodes, so $\tau_\w(T_\ell) = 1 + 5\w(1) \ge 6$.
    Moreover, $\targetbd_\ell > 4\cbd$, as otherwise $\PureBB{\targetbd_\ell}$ has at most $5$ nodes and time at most $5 < 6 \le \tau_\w(T_\ell)$, contradicting minimality.
    Hence, the right child of $u$ must prove $\targetbd_\ell - 2\cbd \in (2\cbd, 3\cbd]$, so $T_\rbd$, which has no cut nodes, is $\PureBB{4\cbd}$ (with $5$ nodes).
    The left child of the root of $T_\rbd$ has remaining bound in $(0, \cbd]$ and two leaf children; replacing these three nodes by a cut node followed by a leaf yields a valid tree, so the minimality of $T_\rbd$ gives $3 \leq 1+\w(1)$, i.e., $\w(1) \ge 2$.
    On the other hand, $\Delta = 1 + \tau_\w(T_\rbd) - \w(1)\big(p(4\cbd) - p(0)\big) = 6 - 4\w(1) \ge 0$ gives $\w(1) \le 3/2$, which is a contradiction.

    \mbox{}\\[-1ex]\noindent
    \emph{Case 2: $\targetbd_G = 2\cbd$.}
        Here, we cannot define $T'$ via the ``two-layers-down'' operation, as $T_G = \PureBB{2\cbd}$ has only one branch node.
    
    \mbox{}\\[-1ex]\noindent
    \emph{Case 2a: $T_\ell$ contains a leaf cut.}
    $T_\ell$ consists of the $\kappa$ cut nodes followed by a branch node with two children (a leaf cut and a leaf).
    Then $\tau_\w(T_\ell) = \sum_{i=0}^{\kappa-1}\w(i) + 3\w(\kappa) +\w(\kappa+1)$. 
    Replacing $T_\ell$ by $\kappa - 1$ cut nodes followed by $\PureBB{4\cbd}$ (with $5$ nodes) yields a subtree that proves the same bound,
    has tree time $\sum_{i=0}^{\kappa-2}\w(i) + 5\w(\kappa-1) \le \tau_\w(T_\ell)$ since $\w$ is nondecreasing, and has two fewer cut nodes, contradicting cut-minimality of $T$. \Cref{fig:exceptional} shows case $\kappa = 1$.

\ifspringer
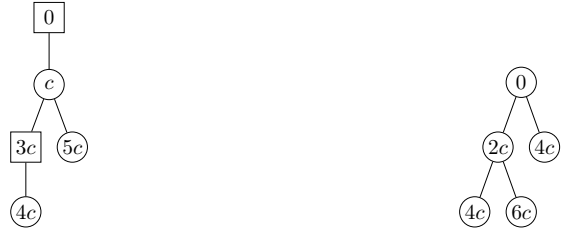
\begin{figure}[t]
\else
\begin{figure}[ht!]
\fi
    \centering
    \subfigure{%
        \scalebox{0.8}{
        \begin{forest}
            for tree = {big node}
            [0,big cut 
                [$c$ 
                    [3$c$, big cut[4$c$]]
                    [5$c$]
                ]
            ]
        \end{forest}
        }
    }
    \hfil
    \subfigure{%
        \scalebox{0.8}{
        \begin{forest}
            for tree = {big node}
            [$0$ 
                [$2c$ [$4c$][$6c$]
                ]
                [
                $4c$
                ]
            ]
        \end{forest}
        }
    }
    
    \caption{Exceptional case when $\targetbd_G=2\cbd$ and $T_\ell$ contains leaf cuts. The tree on the left panel is dominated by the tree on the right, with fewer cut nodes.}
    \label{fig:exceptional}
\end{figure}

\ifspringer
    \mbox{}\\[-1ex]
\else
\fi
    \noindent
    \emph{Case 2b: $T_\ell$ contains no leaf cuts.}
    With no leaf cuts in $T_\ell$, 
    $\tau_\w(T_\ell) = \sum_{i=0}^{\kappa-1}\w(i)+3\w(\kappa)$. Analogous to the ``two-layers-down'' operation, we now employ a ``one-layer-down'' transformation, and we continue to denote the resulting tree by $T'$: the left child of $u$ becomes a branch node whose right child is a copy of $T_1$ and whose left child is the path of $\kappa$ cut nodes followed by a leaf, 
    which suffices as $T_\ell$ proves $\kappa\cbd + 2\cbd$.
    We also continue to use $T''$ for the tree consisting of the $\kappa$ cut nodes followed by $\PureBB{4\cbd}$ (with $5$ nodes), which also proves $\targetbd_u$, as in Case~1.
    Visualizations of these trees are in \Cref{fig:124-base}. 
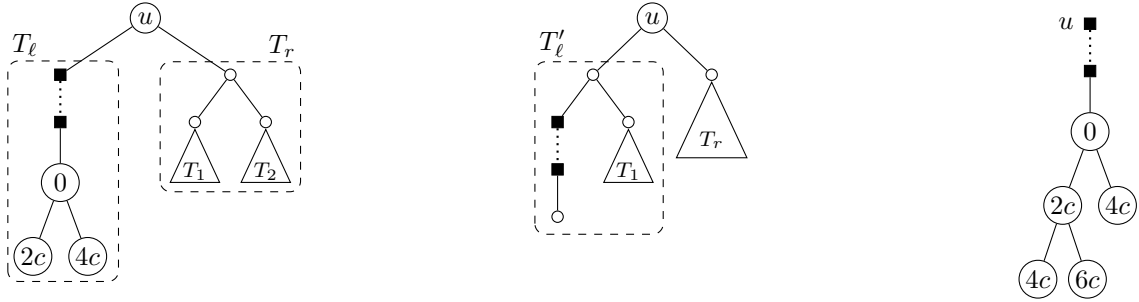
\begin{figure}
  \centering
  \begin{minipage}[t]{0.3\textwidth}\centering
    \begin{forest} baseline, for tree={small node, l=8pt, s sep=14pt}
    [$u$, fill=none, inner sep=1.5pt, s sep=46pt
      [, cut, name=TLroot
        [, cut, cutpath
          [$0$, big node, fill=none, s sep=6pt
            [$2\cbd$, big node, fill=none, name=TLa]
            [$4\cbd$, big node, fill=none, name=TLb]]]]
      [, fill=none, name=TRroot, s sep=22pt
        [, fill=none, tri={T_1}{7mm}, name=TRl]
        [, fill=none, tri={T_2}{7mm}, name=TRr]]]
    \draw[dashed, rounded corners] ([xshift=-4.6mm,yshift=-8.4mm]TRl.south) rectangle ([xshift=4.6mm,yshift=8.8mm]TRr.south);
    \node[anchor=south east, font=\small, inner sep=1.5pt] at ([xshift=4.6mm,yshift=8.8mm]TRr.south) {$T_\rbd$};
    \coordinate (TLBL) at ([xshift=-1.5mm,yshift=-1.5mm]TLa.south west);
    \coordinate (TLR) at ([xshift=1.5mm]TLb.east);
    \coordinate (TLT) at ([yshift=1mm]TLroot.north);
    \draw[dashed, rounded corners] (TLBL) rectangle (TLR |- TLT);
    \node[anchor=south west, font=\small, inner sep=1.5pt] at (TLBL |- TLT) {$T_\lbd$};
    \end{forest}
  \end{minipage}\hfill
  \begin{minipage}[t]{0.38\textwidth}\centering
    \begin{forest} baseline, for tree={small node, l=8pt, s sep=16pt}
    [$u$, fill=none, inner sep=1.5pt, s sep=40pt
      [, fill=none, s sep=22pt
        [, cut [, cut, cutpath [, fill=none, name=TLpleaf]]]
        [, fill=none, tri={T_1}{7mm}, name=TLp1]]
      [, fill=none, tri={T_\rbd}{10mm}]]
    \coordinate (TLpBL) at ([xshift=-3mm,yshift=-1.5mm]TLpleaf.south);
    \coordinate (TLpTR) at ([xshift=4.6mm,yshift=8.8mm]TLp1.south);
    \draw[dashed, rounded corners] (TLpBL) rectangle (TLpTR);
    \node[anchor=south west, font=\small, inner sep=1.5pt] at (TLpBL |- TLpTR) {$T'_\lbd$};
    \end{forest}
  \end{minipage}\hfill
  \begin{minipage}[t]{0.22\textwidth}\centering
    \begin{forest} baseline, for tree={small node, l=8pt}
    [, cut, label={[font=\small]180:{$u$}}
      [, cut, cutpath
        [$0$, big node, fill=none, s sep=6pt
          [$2\cbd$, big node, fill=none, s sep=4pt
            [$4\cbd$, big node, fill=none]
            [$6\cbd$, big node, fill=none]]
          [$4\cbd$, big node, fill=none]]]]
    \end{forest}
  \end{minipage}
  \caption{Case $\targetbd_G = 2\cbd$ with no leaf cuts.
    Left panel: tree $T$.
    Middle panel: tree $T'$, obtained by moving the $\kappa$ cut nodes down a layer.
    Right panel: tree $T''$, obtained by moving the $\kappa$ cut nodes up a layer to $u$.
    Node labels are bounds relative to the node after the cut path.}
  \label{fig:124-base}
\end{figure}
    Denoting by $T'_\ell$ the subtree of $T'$ rooted at the left child of $u$, the minimality of $T$ gives
    \begin{equation*}
        \tau_\w(T'_\ell) = 1 + \tau_{\w}(T_1) + \sum_{i=0}^{\kappa}\w(i) \geq \sum_{i=0}^{\kappa-1}\w(i)+3\w(\kappa) = \tau_\w(T_\ell),
    \end{equation*}
    which implies $1 + \tau_{\w}(T_1) \geq 2\w(\kappa)$.
    Comparing the tree times of $T$ and $T''$, 
    and using $\tau_{\w}(T_r) = 1 + \tau_{\w}(T_1) + \tau_{\w}(T_2) > \tau_{\w}(T_1)$,
    we obtain
    \begin{align*}
        \tau_{\w}(T) - \tau_{\w}(T'') 
        &=\left(1 + \tau_{\w}(T_r) + \sum_{i=0}^{\kappa-1}\w(i)+3\w(\kappa)\right) - \left(\sum_{i=0}^{\kappa-1}\w(i)+5\w(\kappa)\right) \nonumber\\
        &= (1 + \tau_{\w}(T_r)) - 2\w(\kappa) 
        > (1+\tau_{\w}(T_1)) - 2\w(\kappa) \ge 0, %
    \end{align*}
    contradicting the minimality of $T$.

    As every case leads to a contradiction, this completes the proof.
\end{proof}

\Cref{thm:124} demonstrates that the importance of root-cut structure extends beyond symmetric branching. For this special asymmetric setting, not using internal cuts remains optimal uniformly over arbitrary target bounds and time functions. 
What, in particular, is special about this setting? While we do not have a precise answer, computationally, we could not find any instances requiring internal cuts when $\rbd = 2\lbd$, but the current proof does not directly extend to this more general case.

The results derived for the single cut-family, single branching-family model may not easily extend to settings with multiple branching variables or cut families.
For instance, while internal cuts are unnecessary under the individual settings $(\cbd_1,\lbd_1,\rbd_1) = (1,2,4)$ (by \Cref{thm:124}) and $(\cbd_2,\lbd_2,\rbd_2) = (1,3,3)$ (by \Cref{thm:rootcuts}), the existence of an internal-cuts-free minimal-time tree is no longer guaranteed when the two settings are combined, as \Cref{ex:combined-settings} shows.

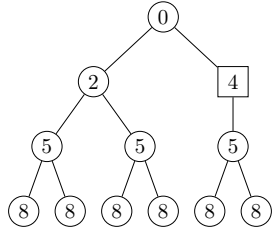
\begin{figure}[ht]
    \centering
    \scalebox{0.8}{%
    \begin{forest}
        for tree = {big node}
        [0
            [2
                [5 [8] [8]]
                [5 [8] [8]]
            ]
            [4, big cut
                [5 [8] [8]]
            ]
        ]
    \end{forest}%
    }
    \caption{The unique minimal-time tree in \Cref{ex:combined-settings}, for $\cbd = 1$, branching improvements $(2,4)$ or $(3,3)$ at each branch node, $\targetbd = 8$, and $\w = (1, a, M, \ldots)$ with $a \in (1, 4/3)$. The root uses $(2,4)$ and all other branch nodes use $(3,3)$; the cut node with bound $4$ is an internal cut. Node labels are bounds, and square nodes are cut nodes.}
    \label{fig:combined-settings}
\end{figure}

\begin{example}\label{ex:combined-settings}
    Consider a single cut family with constant strength $\cbd = 1$ and two branching variables, with bound improvements $(2,4)$ and $(3,3)$, so that each branch node may use either pair.
    Let $\targetbd = 8$ and $\w = (1, a, M, M, \ldots)$, with $a \in (1, 4/3)$ and $M$ sufficiently large, so that at most one cut is used on any root-to-leaf path.
    The tree in \Cref{fig:combined-settings} branches at the root using $(2,4)$; its left child is followed by three branch nodes using $(3,3)$ and four leaves, and its right child is an internal cut followed by a branch node using $(3,3)$ with two leaves.
    Its tree time is $9 + 3a$, and it is the unique minimal-time tree, whereas the best tree without internal cuts has tree time $\min\{1 + 11a,\, 13\}$, attained by a root cut followed by pure branching or by pure branching.
\end{example}

\section{How Close Are Root-Cut Trees to Optimal?}
\label{sec:rootcut-quality}

\Cref{root} studies settings in which some minimal-time tree has no internal cuts,
but generally, minimal-time trees may need local cuts, possibly at any depth (\cref{sec:complex}).
In this section, we quantify how much is lost by restricting cuts to the root:
\cref{sec:notbad} bounds the ratio between the tree time of the best root-cut tree and that of a minimal-time tree for any fixed target bound, and \cref{sec:asymptotic} shows that, as the target bound grows, the best root-cut tree grows at the same rate as a minimal-time tree
and, for constant cuts, is asymptotically optimal.

Throughout this section, we consider a single cut family setting $(\cbd,\lbd,\rbd) \in \cutfunset \times \posreals \times \posreals$ and a generalized time function $\W$ with simple time function $\w$ and separation time function $\septime$ (of the cut family).
Thus, a branch or leaf node with $h$ cut nodes above it has time $\w(h)$, and a cut node with $h$ cut nodes above it has time $\w(h) + \septime(h)$.
Let $C_k \defeq \sum_{h=1}^{k}\cbd(h)$ denote the cumulative bound improvement of $k$ cuts.
For a target bound $\targetbd$, let $\tau^*_{\W}(\targetbd)$ denote the tree time of a minimal-time tree proving $\targetbd$ with respect to $\W$, with no restriction on where cuts are placed, and let $t^*(\targetbd)$ denote the size of a minimal-size tree proving $\targetbd$.
We use $p(z) = \cardinality{\PureBB{z}}$ for $z \in \R$, as in \cref{sec:124}.

For $k \in \nonnegints$, let $\pi(k)$ be the tree time of a path of $k$ cut nodes then a leaf:
\begin{equation*}
    \pi(k) \defeq \sum_{i \in [0,k-1]}\big(\w(i) + \septime(i)\big) + \w(k).
\end{equation*}
A root-cut tree with $k$ cut nodes proves $\targetbd$ if and only if its branching component proves $\targetbd - C_k$.
All of its branch and leaf nodes have time $\w(k)$, so the minimal-time such tree uses $\PureBB{\targetbd - C_k}$ and %
\begin{equation}\label{eq:rootcut-time}
    \tau_{\mathrm{root}}(k) \defeq \sum_{\smash{i=0}}^{\smash{k-1}}\big(\w(i) + \septime(i)\big) + \w(k)\, p(\targetbd - C_k)
    = \pi(k) + \w(k)\big(p(\targetbd - C_k) - 1\big).
\end{equation}
The minimal time of a root-cut tree proving $\targetbd$ is $\tau^*_{\mathrm{root}}(\targetbd) \defeq \min_{k \in \nonnegints} \tau_{\mathrm{root}}(k)$, and we denote by $k_{\mathrm{root}}(\targetbd)$ a minimizer.
We refer to a root-cut tree attaining $\tau^*_{\mathrm{root}}(\targetbd)$ as a \emph{best root-cut tree}.

\subsection{Root-Cut Approximation for a Fixed Target Bound}
\label{sec:notbad}

We first bound the tree time of any tree from below, in terms of the maximum number of cut nodes along a root-to-leaf path.

\begin{lemma}\label{lem:inter}
    Let $T$ be any tree proving $\targetbd$, and let $k$ be the maximum number of cut nodes on a root-to-leaf path of $T$.
    Then $\cardinality{T} \ge k + p(\targetbd - C_k)$ and
    \begin{equation*}
        \treetime_{\W}(T) \ge \pi(k) + p(\targetbd - C_k) - 1.
    \end{equation*}
\end{lemma}
 
\begin{proof}
    Deleting the cut nodes of $T$ leaves a pure branching tree $B$.
    Each leaf of $B$ has at most $k$ cut nodes above it in $T$, which improve its bound by at most $C_k$, so $B$ proves $\targetbd - C_k$ and $\cardinality{B} \ge p(\targetbd - C_k)$.
    As $T$ has at least $k$ cut nodes, $\cardinality{T} \ge k + \cardinality{B}$, which proves the first claim.
    For the second, a root-to-leaf path with $k$ cut nodes has time $\pi(k)$, and each of the remaining $\cardinality{T} - k - 1$ nodes has time at least $1$.
\end{proof}

\Cref{lem:inter} uses the argument of \cite[Lemma~7]{KazLeBSan24_abstractbc-mpb} and \cite[Lemma~4.1]{Jiang26_sharp-logarithmic-thresholds_arxiv}, which move all cut nodes to the root without increasing the size of a tree.
With this, we bound the root-cut approximation for a minimal-time tree.
For a minimal-time tree $T^*$ proving $\targetbd$, let $k^*$ denote the maximum number of cut nodes on a root-to-leaf path of $T^*$, 
and let $B^* \defeq p(\targetbd - C_{k^*}) - 1$.
 
\begin{proposition}\label{prop:notbad}
    For every target bound $\targetbd$,
    \begin{equation*}
        \pi(k^*) + B^* \le \tau^*_{\W}(\targetbd) \le \tau^*_{\mathrm{root}}(\targetbd) \le \min\{p(\targetbd),\; \pi(k^*) + \w(k^*)\,B^*\}.
    \end{equation*}
    Consequently,
    \begin{equation*}
        \frac{\tau^*_{\mathrm{root}}(\targetbd)}{\tau^*_{\W}(\targetbd)}
        \le \min\left\{ \frac{p(\targetbd)}{\pi(k^*) + B^*},\; \frac{\pi(k^*) + \w(k^*)\,B^*}{\pi(k^*) + B^*} \right\}
        \le \w(k^*).
    \end{equation*}
\end{proposition}
 
\begin{proof}
    The first inequality is \cref{lem:inter} applied to $T^*$, and the second holds as root-cut trees are a subset of all trees.
    For the third, \eqref{eq:rootcut-time} gives $\tau_{\mathrm{root}}(0) = p(\targetbd)$, as $\pi(0) = \w(0) = 1$, and $\tau_{\mathrm{root}}(k^*) = \pi(k^*) + \w(k^*)\,B^*$.
    The last inequality follows from $\pi(k^*) \le \w(k^*)\,\pi(k^*)$.
\end{proof}

The factor $\w(k^*)$ is the bound of \cite[Theorem~11]{HanKaz26_abstractbc-ipco}.
\Cref{prop:notbad} improves it to $(\pi(k^*) + \w(k^*)B^*)/(\pi(k^*) + B^*)$ and extends it to worsening cut functions and separation time.
The bound is exact when $\w \equiv 1$ and $\septime \equiv 0$, and when $B^* = 0$, in which case the pure cutting tree with $k^*$ cut nodes is a minimal-time tree.
For constant cuts and converging separation time, the bound tends to $1$ as $\targetbd$ grows (\cref{thm:largeZallcuts}).
We revisit this factor via simulation in \cref{sec:simulation}.

\subsection{Asymptotic Optimality of Root-Cut Trees}\label{sec:asymptotic}

\Cref{prop:notbad} bounds how far the best root-cut tree can be from a minimal-time tree for a fixed target bound.
We now ask how this comparison behaves as the target bound grows.
\citet*{LeBNem17} study the asymptotic growth of the size of a minimal \bb tree, as the target bound $\targetbd$ tends to infinity.
\citet*{Jiang26_sharp-logarithmic-thresholds_arxiv} extends this analysis to the abstract \bc model, focusing on size and the proportion of the bound proved by cut nodes.
We further generalize these results from tree size to tree time, including the time to generate cuts, and show that as the target bound increases, the best root-cut tree grows at the same rate as a minimal-time tree (\cref{lem:growthrate}) and becomes asymptotically optimal for constant cuts (\cref{thm:largeZallcuts}).

Throughout this subsection, we assume that $\w$ and $r$ are both bounded from above.
The sequence $\{\w(i)\}$ is then convergent, as it is nondecreasing and upper bounded.
Define $W_{\max} \coloneqq \lim_{i\rightarrow\infty}\w(i) \geq 1$ and $\sepmax \coloneqq \sup_{h \in \nonnegints} \septime(h) \ge 0$.
Let $\mu^* > 0$ be the unique solution of $e^{-\mu^*\lbd} + e^{-\mu^*\rbd} = 1$, which governs the exponential growth rate of pure branching trees~\cite{LeBNem17,Jiang26_sharp-logarithmic-thresholds_arxiv}.
\citet[Theorem~4.3]{Jiang26_sharp-logarithmic-thresholds_arxiv} shows that minimal-size trees grow at a rate determined by how $C_k$ grows relative to $\log k$.
\Cref{lem:growthrate} shows that, under bounded node and separation times, the same holds for the best root-cut tree and for a minimal-time tree.
A technical point is that we allow the limit $\gamma$ in \cref{lem:growthrate} to be $\infty$, i.e., $C_k/\log k \to \infty$.

\begin{lemma}\label{lem:growthrate}
    For every target bound $\targetbd$,
    \begin{equation*}
        t^*(\targetbd) \le \tau^*_{\W}(\targetbd) \le \tau^*_{\mathrm{root}}(\targetbd) \le (W_{\max} + \sepmax)\, t^*(\targetbd).
    \end{equation*}
    Further, if the limit $\gamma \defeq \lim_{k\rightarrow\infty} C_k/\log k$ exists in $[0,\infty]$, then
    \begin{equation*}
        \lim_{\targetbd\rightarrow\infty} \frac{\log(\tau^*_{\mathrm{root}}(\targetbd))}{\targetbd}
        =
        \lim_{\targetbd\rightarrow\infty} \frac{\log(\tau^*_{\W}(\targetbd))}{\targetbd} 
        = 
        \lim_{\targetbd\rightarrow\infty}\frac{\log(t^*(\targetbd))}{\targetbd}
        =
        \frac{\mu^*}{1+\gamma\mu^*},
    \end{equation*}
    where the right-hand side of the last equality is taken as $0$ when $\gamma = \infty$.
\end{lemma}

\begin{proof}
    Since every node time is at least $1$, the tree time of any tree is at least its size, so $t^*(\targetbd) \le \tau^*_{\W}(\targetbd)$, and $\tau^*_{\W}(\targetbd) \le \tau^*_{\mathrm{root}}(\targetbd)$ holds as root-cut trees are a subset of all trees.
    Let $T$ be a minimal-size tree proving $\targetbd$, and let $k$ be the maximum number of cut nodes on a root-to-leaf path of $T$.
    By \cref{lem:inter}, the root-cut tree consisting of $k$ cut nodes followed by $\PureBB{\targetbd - C_k}$, which proves $\targetbd$, has size at most $\cardinality{T}$, so it is also a minimal-size tree (as in \cite[Lemma~7]{KazLeBSan24_abstractbc-mpb} and \cite[Lemma~4.1]{Jiang26_sharp-logarithmic-thresholds_arxiv}).
    As every node time is at most $W_{\max} + \sepmax$, its tree time is at most $(W_{\max} + \sepmax)\, t^*(\targetbd)$, which proves the last inequality.
    When $\gamma$ exists, the limit of $\log(t^*(\targetbd))/\targetbd$ exists and equals $\mu^*/(1+\gamma\mu^*)$ by \cite[Theorem~4.3]{Jiang26_sharp-logarithmic-thresholds_arxiv}. 
    Taking logarithms in the first claim, dividing by $\targetbd$, and using $\log(W_{\max} + \sepmax)/\targetbd \to 0$ completes the proof.
\end{proof}

Moreover, by \cref{prop:notbad}, the best root-cut tree is always within a constant factor of a minimal-time tree:
$\tau^*_{\mathrm{root}}(\targetbd) \le \w(k^*)\,\tau^*_{\W}(\targetbd) \le W_{\max}\,\tau^*_{\W}(\targetbd)$, even when the separation time is unbounded.
When $\gamma$ exists, the best root-cut tree, a minimal-time tree, and a minimal-size tree have the same exponential growth rate in the target bound.
This lemma is particularly informative when tree size grows exponentially, i.e., when $\gamma < \infty$, 
such as for harmonically worsening cuts $\cbd(h) = \cbd_0/h$ ($\gamma = \cbd_0$) or when $\sum_{h \ge 1} \cbd(h) < \infty$ ($\gamma = 0$).
In contrast, tree size grows linearly under a constant cut function $\cbd > 0$, and the minimal tree size satisfies $t^*(\targetbd) \le \lceil \targetbd/\cbd\rceil + 1$, so the limits in \cref{lem:growthrate} equal $0$ (this is captured by $\gamma = \infty$).
In this case, 
\cref{thm:largeZallcuts} gives a sharper result:
the best root-cut tree is asymptotically optimal, and its root cuts prove asymptotically all of the target bound.

The key ingredient is the pure cutting tree, the root-cut tree with no branch nodes.
For a constant cut function $\cbd > 0$ and a target bound $\targetbd$, let $K_{\targetbd} \defeq \lceil \targetbd/\cbd \rceil$, let $T_{\mathrm{cut}}(\targetbd)$ denote the pure cutting tree proving $\targetbd$, consisting of $K_{\targetbd}$ cut nodes followed by a leaf, and let
\begin{equation*}
    \tau_{\mathrm{cut}}(\targetbd) \defeq \treetime_{\W}\big(T_{\mathrm{cut}}(\targetbd)\big) = \pi(K_{\targetbd})
\end{equation*}
denote its tree time, which we compute asymptotically in \cref{lem:limit:w/c}.
For \cref{lem:limit:w/c,thm:largeZallcuts}, we also assume that the separation time converges, $\septime(h) \to \sepinf$ as $h \to \infty$. 
This holds, e.g., for a constant separation time $\septime \equiv \sigma$, with $\sepinf = \sigma$, or for no separation time, with $\sepinf = 0$.

\begin{lemma}\label{lem:limit:w/c}
    Suppose that $\septime(h) \to \sepinf$ as $h \to \infty$.
    Then $\pi$ is strictly increasing and $\lim_{k\to\infty}\pi(k)/k = W_{\max} + \sepinf$.
    Hence, for a constant cut function $\cbd > 0$,
    \begin{equation*}
        \lim_{\targetbd\rightarrow\infty}\frac{\tau_{\mathrm{cut}}(\targetbd)}{\targetbd} = \frac{W_{\max} + \sepinf}{\cbd}.
    \end{equation*}
\end{lemma}
 
\begin{proof}
    For every $k$, $\pi(k+1) - \pi(k) = \septime(k) + \w(k+1) \ge 1$.
    Since $\w(i) + \septime(i) \to W_{\max} + \sepinf$, the averages $\frac{1}{k}\sum_{i=0}^{k-1}\big(\w(i) + \septime(i)\big)$ converge to the same limit, and $\w(k)/k \le W_{\max}/k \to 0$, so $\pi(k)/k \to W_{\max} + \sepinf$.
    The last claim follows from $\tau_{\mathrm{cut}}(\targetbd) = \pi(K_\targetbd)$ and $K_\targetbd/\targetbd \to 1/\cbd$.
\end{proof}

\begin{theorem}\label{thm:largeZallcuts}
    Let $\cbd > 0$ be a constant cut function, and suppose that $\septime(h) \to \sepinf$ as $h \to \infty$. Then
    \begin{equation*}
        \lim_{\targetbd\rightarrow\infty}\frac{\tau^*_{\mathrm{root}}(\targetbd)}{\targetbd} = \lim_{\targetbd\rightarrow\infty}\frac{\tau^*_{\W}(\targetbd)}{\targetbd} = \frac{W_{\max} + \sepinf}{\cbd},
        \quad\text{so}\quad
        \lim_{\targetbd\rightarrow\infty}\frac{\tau^*_{\mathrm{root}}(\targetbd)}{\tau^*_{\W}(\targetbd)} = 1.
    \end{equation*}
    Moreover, $\targetbd - \cbd\, k_{\mathrm{root}}(\targetbd) \in O(\log \targetbd)$ and, for every choice of minimal-time tree $T^*$ proving $\targetbd$, $\targetbd - \cbd\, k^* \in O(\log \targetbd)$.
    In particular, the root cuts of a best root-cut tree, and the cuts on some root-to-leaf path of a minimal-time tree, prove a fraction of the target bound that tends to $1$.
\end{theorem}
\begin{proof}
    Let $T^*$ be any minimal-time tree proving $\targetbd$, with $k^*$ the maximum number of cut nodes on a root-to-leaf path of $T^*$.
    Applying \cref{lem:inter} to $T^*$, and since $T_{\mathrm{cut}}(\targetbd)$ is a root-cut tree,
    \begin{equation}\label{eq:largeZ-sandwich}
        \pi(k^*) + B^*
        =
        \pi(k^*) + p(\targetbd - \cbd k^*) - 1 \le \tau^*_{\W}(\targetbd) \le \tau^*_{\mathrm{root}}(\targetbd) \le \pi(K_{\targetbd}).
    \end{equation}
    As $\pi$ is strictly increasing (\cref{lem:limit:w/c}) and $p(z) \ge 1$ for all $z$, this implies $k^* \le K_\targetbd$.
    Every node of $\PureBB{z}$ at depth less than $z/\rbd$ has bound less than $z$, so $\PureBB{z}$ contains a complete binary tree of depth $\lceil z/\rbd \rceil$, and $p(z) \ge 2^{z/\rbd}$.
    Hence, by \eqref{eq:largeZ-sandwich} and $\pi(K_\targetbd) \le (W_{\max} + \sepmax)K_\targetbd + W_{\max}$, we have $\targetbd - \cbd k^* \le \rbd \log_2 \pi(K_\targetbd) \in O(\log \targetbd)$.
    Similarly, a best root-cut tree has time $\pi(k_{\mathrm{root}}) + \w(k_{\mathrm{root}})\big(p(\targetbd - \cbd k_{\mathrm{root}}) - 1\big) \le \pi(K_\targetbd)$, so $k_{\mathrm{root}} \le K_\targetbd$ (as $\w \ge 1$) and $\targetbd - \cbd k_{\mathrm{root}} \in O(\log \targetbd)$.
 
       Therefore, $\targetbd - O(\log \targetbd) \le \cbd k^* \le \cbd K_\targetbd < \targetbd + \cbd$, so, as $\targetbd \to \infty$, we have $k^* \to \infty$ (from the lower bound) and $\cbd k^*/\targetbd \to 1$.
        Since $k^* \to \infty$, \cref{lem:limit:w/c} gives $\pi(k^*)/k^* \to W_{\max} + \sepinf$, so, as $\targetbd \to \infty$,
    \begin{equation*}
        \frac{\pi(k^*)}{\targetbd} 
        = \frac{\pi(k^*)}{k^*} \cdot \frac{k^*}{\targetbd} 
        \to \frac{W_{\max} + \sepinf}{\cbd},
    \end{equation*}
    and likewise $\pi(K_\targetbd)/\targetbd \to (W_{\max} + \sepinf)/\cbd$.
    Dividing \eqref{eq:largeZ-sandwich} by $\targetbd$ proves both limits, and hence that their ratio tends to $1$.
\end{proof}

\Cref{thm:largeZallcuts} shows that, for constant cuts, restricting cuts to the root is asymptotically optimal and the root cuts prove all but a logarithmic portion of the target bound. 
For $\w \equiv 1$, no separation time, and constant cuts ($\gamma = \infty$), \cref{thm:largeZallcuts} recovers \cite[Corollary~4.4]{Jiang26_sharp-logarithmic-thresholds_arxiv}, which states that the fraction of the target bound proved by the root cuts of a minimal-size root-cut tree tends to $1$, and extends it to every minimal-size tree, in that the cuts on some root-to-leaf path prove this fraction.

The proof of \cref{thm:largeZallcuts} implies that the factor $(\pi(k^*) + \w(k^*)B^*)/(\pi(k^*) + B^*)$ from \cref{prop:notbad} tends to $1$,%
\footnote{By \eqref{eq:largeZ-sandwich}, $B^* \le \pi(K_\targetbd) - \pi(k^*) 
\le (W_{\max}+\sepmax)(K_\targetbd - k^*) + W_{\max}$,
which is in $O(\log \targetbd)$, while $\pi(k^*)$ grows linearly.}
whereas the factor $\w(k^*)$ of \cite[Theorem~11]{HanKaz26_abstractbc-ipco} tends to $W_{\max}$.
We emphasize that using node times rather than tree sizes alone is what makes the bound in \cref{prop:notbad} asymptotically tight.
 Although larger separation time may cause fewer cut nodes for a fixed target bound (see \cref{sec:separation}), \cref{thm:largeZallcuts} shows that root cuts still asymptotically suffice for constant cuts.

\section{Internal Cuts Can Be Necessary}
\label{local}
In this section, we explore general asymmetric branching settings with a single constant-strength cut family in which internal cuts can be necessary in any minimal-time tree.
This is the common practice in modern solvers, where cuts generated beyond the root node play a critical role in reducing solving time and improving performance~\cite{BerFraHen25_learning-to-use-local-cuts}.

\subsection{Minimal-Time Trees Can Be Complex}\label{sec:complex}
Unlike the symmetric branching setting or the special case analyzed in \Cref{sec:124}, where a minimal-time tree exists with cuts only in one of two locations (at the root or at a leaf), general asymmetric branching settings together with special choices of time functions can yield far more complex optimal structures, making prescribing best cut locations nontrivial. We construct an intuitive example in \Cref{lem-complex} to illustrate the complexity of minimal-time trees under asymmetric branching settings.

\begin{lemma}\label{lem-complex}
    Let $(\cbd,\lbd,\rbd) = (2,3,4)$. %
    For any depth $d \in \nonnegints$,
    there exists a simple time function $\w$
    such that any minimal-time tree proving bound $\targetbound \defeq 3 d + 5$ will have a cut node at depth $d$ along the leftmost root-to-leaf path.
\end{lemma}
\begin{proof}
    For $d \in \nonnegints$, let $\bar\w_d \defeq 2 - 1/p(3d+9)$, pick any $\w(1) \in [\bar\w_d,2)$, and set $\w(k) = M$ for all $k \ge 2$ with $M$ sufficiently large that every minimal-time tree has at most one cut node on each root-to-leaf path.
    We proceed by induction on $d$.
    For $d = 0$ and $\targetbd = 5$, a root cut followed by a branch node with two leaves has tree time $1 + 3\w(1)$.
    As $\w(1) < 2$, this is less than the time 7 of $\PureBB{5}$, the time $3 + 2\w(1)$ of a branch node whose children are both leaf cuts, and the time $5 + \w(1)$ of a branch node with one leaf-cut child, so every minimal-time tree has a cut node at its root.
 
    Now suppose that the lemma holds for some $d \ge 0$ and every $\w(1) \in [\bar\w_d, 2)$, and consider $\targetbd = 3d+8$ and $\w(1) \in [\bar\w_{d+1}, 2)$.
    In a tree whose root is a cut node, the remaining nodes form a pure branching tree proving $\targetbd - 2$, each with time $\w(1)$, so its tree time is at least $1 + \w(1)\,p(\targetbd - 2) > 2p(\targetbd - 2)$, as $\w(1) \ge \bar\w_{d+1} > 2 - 1/p(\targetbd - 2)$.
    On the other hand, $p(\targetbd) = 1 + p(\targetbd - 3) + p(\targetbd - 4) \le 2p(\targetbd - 2)$, as $p$ is nondecreasing and $p(\targetbd - 4) < p(\targetbd - 2)$ (the leaf of $\PureBB{3d+4}$ with bound $3d+4$ is a branch node of $\PureBB{3d+6}$).
    Hence, $\PureBB{\targetbd}$ is faster than every tree with a root cut, so the root of every minimal-time tree is a branch node.
    By \cite[Lemma~24]{KazLeBSan24_abstractbc-mpb}, the left subtree of the root is a minimal-time tree proving $3d+5$, and, as $\w(1) \ge \bar\w_{d+1} \ge \bar\w_d$, the induction hypothesis gives a cut node at depth $d$ on its leftmost path, i.e., at depth $d+1$ on the leftmost path of the whole tree.
\end{proof}

In the asymmetric setting of \Cref{lem-complex}, minimal-time trees have arbitrarily deep internal cuts. However, in the construction, the internal cut is always two layers above the bottom of the tree. In other words, while the depth is arbitrary in this setting, the placement of cut nodes is not.

\subsection{Adding At Most One Cut Node Per Root-To-Leaf Path}\label{sec:1cut}

Typically, solvers limit the number of cut rounds, as successive rounds often yield diminishing bound improvements \cite{DeyKazLodMun22} while substantially increasing problem size and solution time. Moreover, excessive cuts can harm numerical stability and lead to ill-conditioned subproblems that degrade solver performance \cite{ZanFisBal11,MilRalSte18_exploring-numerics-branch-and-cut}.
In this subsection, we investigate necessary and sufficient conditions under which some minimal-time tree has an internal cut. We restrict our focus to the \emph{one-cut} setting, defined as a configuration where at most one round of cuts is added in each root-leaf path.
This assumption simplifies the analysis while preserving key structural insights.

Note that if $\w(1) > 2$ under the one-cut setting, no cut nodes would appear in any minimal-time tree,
as branch nodes provide both better bound improvements and smaller processing time when $\w(1)$ is large. 
Meanwhile, when $\w(1) = 1$, there always exists a minimal-time root-cut tree. 
Accordingly, we focus our analysis on $\w(1) \in (1,2)$. 
Under this assumption,
a branch node in a minimal-time tree cannot have two cut node children.

Other than the time function, the relationship among $\targetbound$, $\cbd$, $\lbd$, and $\rbd$ also determines whether cut nodes will help. For instance, if $\cbd=1$ and all other parameters are even numbers, cut nodes would not appear in any minimal-time tree, for any value of $\w(1)$. 
In \Cref{lem:existcuts}, we show that cuts are strictly beneficial if and only
if the pure branching tree has a node with at least 1 and at most $\cbd$ bound left to prove.

\begin{restatable}{lemma}{existcuts}
\label{lem:existcuts}
    Given a one-cut setting $(\targetbound,\cbd,\lbd,\rbd) \in \posints^4$ and $\w(1)\in(1,2)$, every minimal-time tree proving $\targetbd$ has a cut node if and only if there exist $k_1, k_2 \in \nonnegints$ such that $k_1 \lbd + k_2 \rbd \in \{\targetbound-\cbd, \ldots, \targetbound-1\}$.
\end{restatable}
 
\begin{proof}
    Suppose first that such a pair $(k_1, k_2)$ exists, and let $u$ be the node of $\PureBB{\targetbd}$ reached from the root by $k_1$ left and $k_2$ right branches.
    This node exists and is a branch node, as its bound $k_1\lbd + k_2\rbd$ and those of its ancestors are less than $\targetbd$, and its children are leaves, as $k_1\lbd + k_2\rbd + \lbd \ge \targetbd - \cbd + \lbd \ge \targetbd$.
    Replacing $u$ and its two children by a cut node followed by a leaf yields a tree that still proves $\targetbd$, as $k_1\lbd + k_2\rbd + \cbd \ge \targetbd$, and whose tree time is smaller by $3 - (1 + \w(1)) > 0$.
    As $\PureBB{\targetbd}$ is the smallest pure branching tree, no pure branching tree is minimal-time.
 
    Conversely, if no such pair exists, let $v$ be a cut node of a tree $T$ proving $\targetbd$.
    Any leaf below $v$ is reached from the root by some $k_1$ left and $k_2$ right branches, and, as $v$ is the only cut on its path, its bound is $k_1\lbd + k_2\rbd + \cbd \ge \targetbd$.
    Hence, $k_1\lbd + k_2\rbd \ge \targetbd - \cbd$, so $k_1\lbd + k_2\rbd \ge \targetbd$ by assumption, as all parameters are integers.
    Thus, removing $v$ (connecting its child to its parent) yields a strictly faster tree that still proves $\targetbd$, so $T$ is not minimal-time.
\end{proof}

\looseness=-1
When the condition of \Cref{lem:existcuts} holds, two trees serve as extreme cases.
The \emph{one-root-cut tree} is a cut node followed by $\PureBB{\targetbd - \cbd}$, with tree time $1 + p(\targetbd - \cbd)\,\w(1)$.
The \emph{minimal leaf-cut tree} branches as $\PureBB{\targetbd - \cbd}$ and adds a cut above each of its $b_{\text{leaf}}(\targetbd)$ leaves with remaining bound in $[\cbd]$, with tree time $p(\targetbd - \cbd) + b_{\text{leaf}}(\targetbd)\,\w(1)$.
As $b_{\text{leaf}}(\targetbd) < p(\targetbd - \cbd)$, the first is faster for small $\w(1)$ and the second for large $\w(1)$.
We write $b_{\text{leaf}}$ for $b_{\text{leaf}}(\targetbd)$ when the target bound is clear.

We next prove that \Cref{lem:existcuts} constrains the \emph{relative} placement of cuts in sibling subtrees.
Concretely, in a one-cut setting, if one child of a branch-then-cut node begins with an internal cut (a one-root-cut subtree), the other child will also contain a node
where a cut saves time --- so the sibling cannot be cut-free.
\Cref{lem:impossiblesub} formalizes this observation.

\begin{restatable}{lemma}{impossiblesub}
\label{lem:impossiblesub}
    Given a one-cut setting $(\targetbd,\cbd,\lbd,\rbd) \in \posints^4$ and $\w(1) \in(1, 2)$, 
    a minimal-time tree cannot contain a branch-then-cut node 
    where one child is a pure branching subtree.
\end{restatable}

\begin{proof}
    Suppose, for the sake of contradiction, that a minimal-time tree $T$ contains a branch-then-cut node $u$ whose one child is an internal cut node while the other is a branch node followed by a pure branching subtree. Without loss of generality, assume the left branch $T_\ell$ (rooted at the left child of $u$) follows a cut-then-branch framework while the right branch $T_\rbd$ (rooted at the right child of $u$) is a pure branching subtree.
    Let $\bar{\targetbd} \defeq \targetbd - z_T(u) - \lbd$.
    As there is no cut above $u$ in the one-cut setting, \cite[Lemma~24]{KazLeBSan24_abstractbc-mpb} implies that $T_\ell$ and $T_\rbd$ are minimal-time with respect to the same $\w$, for target bounds $\bar{\targetbd}$ and $\bar\targetbd+\lbd-\rbd$, respectively.
    In particular, $T_\ell$ is a cut node followed by $\PureBB{\bar\targetbd - \cbd}$, and $\bar\targetbd > \cbd$ as this cut is internal.
    By \Cref{lem:existcuts}, the existence of the cut node in $T_\ell$ implies that there exists a pair $(k_1,k_2)$ satisfying $k_1 \lbd + k_2 \rbd \in \{\bar{\targetbd}-\cbd,\ldots, \bar{\targetbd}-1\}$.

    We first claim that some such pair has $k_2 \ge 1$.
    Otherwise, as $\cbd \le \lbd$, exactly one node $x$ of $\PureBB{\bar\targetbd}$ has remaining bound in $\{1,\ldots,\cbd\}$, namely the node on its leftmost path with bound $k_1\lbd$, and $x$ is not the root as $\bar\targetbd > \cbd$.
    Thus, $\PureBB{\bar\targetbd - \cbd}$ has exactly one fewer branch node than $\PureBB{\bar\targetbd}$, so $p(\bar\targetbd - \cbd) = p(\bar\targetbd) - 2$, and $\tau_\w(T_\ell) = 1 + \w(1)\big(p(\bar\targetbd) - 2\big)$.
    Replacing the branching at $x$ in $\PureBB{\bar\targetbd}$ by a cut node followed by a leaf instead yields a tree proving $\bar\targetbd$ with tree time $p(\bar\targetbd) - 2 + \w(1)$, which is smaller than $\tau_\w(T_\ell)$ by $(\w(1) - 1)\big(p(\bar\targetbd) - 3\big) > 0$, as $p(\bar\targetbd) \ge 5$, contradicting the minimality of $T_\ell$.

    Thus, some pair has $k_2 \geq 1$, and the node of $T_\rbd$ reached from its root by $k'_1 = k_1 + 1$ left branches and $k'_2 = k_2 - 1$ right branches is a branch node whose children are both leaves, as $k'_1\lbd + k'_2\rbd = k_1\lbd + k_2\rbd + \lbd - \rbd \in\{\bar\targetbd+\lbd-\rbd-\cbd, \dots, \bar\targetbd+\lbd-\rbd-1\}$.
    Replacing this branch node and its two leaves in $T_\rbd$ by a cut node followed by a leaf yields a tree $T'_\rbd$ that proves the same bound, with $\tau_\w(T'_\rbd) - \tau_\w(T_\rbd) = 1 + \w(1) - 3 < 0$ as $\w(1)<2$, contradicting the minimality of $T_\rbd$.
    The case in which $T_\rbd$ follows the cut-then-branch framework and $T_\ell$ is a pure branching tree is analogous, with the roles of left and right branches exchanged (using $\cbd \le \rbd$ and the rightmost path).
\end{proof}

\looseness=-1
Specifically, \Cref{lem:impossiblesub} implies that, if $u$ is a deepest branch-then-cut node in a minimal-time tree in a one-cut setting, then one child subtree starts with a root cut and the other must have leaf cuts.
Both children have remaining bound greater than $\cbd$: the internal cut is followed by a branch node, and if the other child had remaining bound at most $\cbd$, replacing its subtree by a cut node then a leaf would be faster.
Hence, the remaining bound at $u$ is greater than $\rbd + \cbd$.

We now turn our attention to characterizing the structure of a minimal-time tree and, subsequently, the placement of the cut nodes. 
\Cref{lem:pattern} characterizes how the structure of a minimal-time tree evolves as $\w(1)$ grows. 

\begin{restatable}{lemma}{pattern}
\label{lem:pattern}
    In a one-cut setting $(\targetbound, \cbd, \lbd, \rbd) \in \posints^4$ satisfying the condition in \Cref{lem:existcuts}, the values of $\w(1) \in (1,2)$ for which some minimal-time tree is the one-root-cut tree, has an internal cut, or is a leaf-cut tree form three (possibly empty) intervals that appear in this order and overlap only at their endpoints.
    In particular, as $\w(1)$ increases, the minimal-time structure can only change from one-root-cut to internal-cut to leaf-cut, possibly skipping a stage.
\end{restatable}

\begin{proof}
    Let $T$ be a tree that is minimal-time for some $\w(1) \in (1,2)$.
    By \Cref{lem:existcuts}, $T$ has a cut node, and its tree time is $a + b\,\w(1)$, where $b$ is the number of nodes below a cut and $a \ge 1$ is the number of other nodes.
    The time of $T$ is at most the time of the one-root-cut tree, so $b \le p(\targetbd - \cbd)$, with equality only if $a = 1$, i.e., $T$ is the one-root-cut tree.

    In $T$, a node with no cut above it and remaining bound in $\{1,\ldots,\cbd\}$ is a cut node followed by a leaf, as otherwise it is a branch node with two leaves, and replacing these three nodes by a cut node followed by a leaf is faster when $\w(1) < 2$.
    For each of the $b_{\text{leaf}}$ leaves of $\PureBB{\targetbd - \cbd}$ with remaining bound in $\{1,\ldots,\cbd\}$, follow the same branches from the root of $T$, skipping over any cut node; each node passed has remaining bound greater than $\cbd$, or greater than $0$ if a cut was passed, so none is a leaf.
    If a cut was passed, the node reached is a leaf below a cut; otherwise, it is a cut node followed by a leaf.
    This gives $b_{\text{leaf}}$ distinct leaves below cuts, so $b \ge b_{\text{leaf}}$.
    If $T$ has an internal cut, its child is a further node below a cut, so $b > b_{\text{leaf}}$.
    Otherwise, every cut node of $T$ is followed by a leaf, so no node with remaining bound greater than $\cbd$ is a cut node, and $T$ branches exactly at these nodes; thus, $T$ is the minimal leaf-cut tree, with $b = b_{\text{leaf}}$.

    Hence, the one-root-cut tree has slope $p(\targetbd - \cbd)$, the minimal leaf-cut tree has slope $b_{\text{leaf}}$, and every other tree that is minimal-time for some $\w(1) \in (1,2)$ has an internal cut and slope strictly in between.
    If trees with times $a_1 + b_1\w(1)$ and $a_2 + b_2\w(1)$ are minimal-time at $\omega_1 < \omega_2$, respectively, then adding $a_1 + b_1 \omega_1 \le a_2 + b_2 \omega_1$ and $a_2 + b_2 \omega_2 \le a_1 + b_1 \omega_2$ gives $(b_2 - b_1)(\omega_2 - \omega_1) \le 0$, so the slope of a minimal-time tree is nonincreasing in $\w(1)$.
    Thus, the one-root-cut tree is minimal-time exactly on an initial interval of $(1,2)$, the minimal leaf-cut tree exactly on a final interval, and trees with internal cuts only in between.
\end{proof}

Let $\alpha(\targetbd)$ be the value of $\w(1)$ at which the one-root-cut and minimal leaf-cut trees proving $\targetbd$ require the same processing time (for fixed $(\cbd,\lbd,\rbd)$).
Concretely,
    \begin{equation*}
    \alpha(\targetbd) \defeq \dfrac{p(\targetbd - \cbd) - 1}{p(\targetbd - \cbd) - b_{\text{leaf}}(\targetbd)}.
    \end{equation*}
For example, with $(\cbd,\lbd,\rbd) = (2,3,4)$, $\alpha(5) = 2$ in the base case of \Cref{lem-complex}, where the one-root-cut tree is minimal-time for all $\w(1) \in (1,2)$.
\Cref{lem:pattern} suggests that this threshold plays a crucial role in determining the architecture of the minimal-time tree, which \Cref{coro:free} makes precise.

\begin{corollary}\label{coro:free}
    Given a one-cut setting $(\targetbd, \cbd, \lbd, \rbd) \in \posints^4$ with $\alpha(\targetbd) \in (1,2)$, if some minimal-time tree proving $\targetbd$ has an internal cut for some $\w(1) \in (1,2)$, then so does some minimal-time tree proving $\targetbd$ for $\w(1) = \alpha(\targetbd)$.
\end{corollary}
\begin{proof}
    Let $\alpha \defeq \alpha(\targetbd)$.
    Define $R(\omega) \defeq 1 + p(\targetbd - \cbd)\,\omega$ and $L(\omega) \defeq p(\targetbd - \cbd) + b_{\text{leaf}}\,\omega$ for the tree times of the one-root-cut and minimal leaf-cut trees at $\w(1) = \omega$, so that $R(\alpha) = L(\alpha)$, and let $I(\omega) = a + b\,\omega$ be the tree time of an internal-cut tree that is minimal-time at $\w(1) = \omega_0 \in (1,2)$, so $I(\omega_0) \le \min\{R(\omega_0), L(\omega_0)\}$.
    By the proof of \cref{lem:pattern}, $b_{\text{leaf}} < b < p(\targetbd - \cbd)$, so $I - R$ is decreasing and $I - L$ is increasing.
    If $\alpha \ge \omega_0$, then $I(\alpha) - R(\alpha) \le I(\omega_0) - R(\omega_0) \le 0$,
    and otherwise, $I(\alpha) - L(\alpha) \le I(\omega_0) - L(\omega_0) \le 0$.
    In both cases, $I(\alpha) \le R(\alpha) = L(\alpha)$.
    Every one-root-cut tree has time at least $R(\alpha)$, every leaf-cut tree has time at least $L(\alpha)$ (as $\alpha < 2$), and no pure branching tree is minimal-time by \cref{lem:existcuts}, so at $\w(1) = \alpha$, either this internal-cut tree is minimal-time, or a faster tree exists, which then has an internal cut.
\end{proof}

Together, \Cref{lem:pattern,coro:free} give the following picture.
If, for different values of $\w(1) \in (1,2)$, there are minimal-time trees of all three types (one-root-cut, internal-cut, and leaf-cut), then there exist $\omega_1 \le \alpha(\targetbd) \le \omega_2$ such that a minimal-time one-root-cut tree exists for $\w(1) \in (1, \omega_1]$, a minimal-time internal-cut tree exists for $\w(1) \in [\omega_1, \omega_2]$, and a minimal-time leaf-cut tree exists for $\w(1) \in [\omega_2, 2)$.

Finally, by combining \Cref{lem:existcuts,lem:impossiblesub,lem:pattern} and \Cref{coro:free}, we can determine whether some value of $\w(1)$ yields a minimal-time tree with internal cuts.

\begin{theorem}\label{thm:alpha}
    Given a one-cut setting $(\targetbd,\cbd,\lbd,\rbd) \in \posints^4$ with $\alpha(\targetbd) \in (1,2)$, there exists $\w(1) \in (1,2)$ for which some minimal-time tree proving $\targetbd$ has an internal cut if and only if there are $k_1, k_2 \in \nonnegints$ with $\targetbd' \defeq \targetbd - k_1\lbd - k_2\rbd > \rbd + \cbd$ such that
    \begin{equation*}
        \alpha(\targetbd) \le \max\{\alpha(\targetbd' - \lbd), \alpha(\targetbd' - \rbd)\}.
    \end{equation*}
    In this case, some minimal-time tree proving $\targetbd$ has an internal cut for $\w(1) = \alpha(\targetbd)$.
\end{theorem}

\begin{proof}
    For convenience, define $\alpha \defeq \alpha(\targetbd)$ and, for the value of $\targetbd'$ under consideration, $\alpha_\beta \defeq \alpha(\targetbd' - \beta)$ for $\beta \in \{\lbd, \rbd\}$.
    Suppose first that some minimal-time tree proving $\targetbd$ has an internal cut for some $\w(1) \in (1,2)$.
    By \Cref{coro:free}, some minimal-time tree $T$ for $\w(1) = \alpha$ has an internal cut, so $T$ has a branch-then-cut node.
    Let $u$ be a deepest branch-then-cut node of $T$, reached from the root by $k_1$ left and $k_2$ right branches.
    As $u$ has a cut below it, there is no cut above it, so its remaining bound is $\targetbd' = \targetbd - k_1\lbd - k_2\rbd$, and $\targetbd' > \rbd + \cbd$ by the discussion after \Cref{lem:impossiblesub}.
    Let $\beta \in \{\lbd,\rbd\}$ be such that the child of $u$ with remaining bound $\targetbd' - \beta$ is an internal cut node.
    By \cite[Lemma~24]{KazLeBSan24_abstractbc-mpb}, the subtree rooted at this child is a minimal-time tree proving $\targetbd' - \beta$, and, as it is a cut node followed by a pure branching tree, it is the one-root-cut tree proving $\targetbd' - \beta$.
    By minimality, the time of this subtree is at most the time of the minimal leaf-cut tree proving $\targetbd' - \beta$, showing that $\alpha \le \alpha_\beta \le \max\{\alpha_\lbd,\alpha_\rbd\}$.

    Conversely, suppose that such $k_1$ and $k_2$ exist, and let $\w(1) = \alpha$.
    As the one-root-cut and minimal leaf-cut trees then have the same tree time, some minimal-time tree $T$ has a branch node at its root (the minimal leaf-cut tree, if the one-root-cut tree is minimal-time).
    Follow $k_1$ left and $k_2$ right branches from the root of $T$, passing through any cut node, to reach a node $u$.
    Each node on the way has remaining bound at least $\targetbd' > \rbd + \cbd$, or greater than $\rbd$ if a cut was passed, so none is a leaf.
    If one of these nodes, including $u$, is a cut node, then it is not the root, and its child is a branch node, as it has positive remaining bound and is not a second cut on the path; thus, $T$ has an internal cut.
    Otherwise, $u$ is a branch node with no cut above it and remaining bound $\targetbd'$.
    Let $\beta \in \{\lbd,\rbd\}$ with $\alpha \le \alpha_\beta$, and let $T_\beta$ be the subtree rooted at the child of $u$ with remaining bound $\targetbd' - \beta$, which is a minimal-time tree proving $\targetbd' - \beta$ by \cite[Lemma~24]{KazLeBSan24_abstractbc-mpb}.
    As $\alpha_\beta \ge \alpha > 1$, the minimal leaf-cut tree proving $\targetbd' - \beta$ has more than one cut, so the condition of \Cref{lem:existcuts} holds for $\targetbd' - \beta$, and, by \Cref{lem:pattern} and its proof, $T_\beta$ has an internal cut, is the one-root-cut tree, or is the minimal leaf-cut tree proving $\targetbd' - \beta$.
    If $T_\beta$ has an internal cut, so does $T$.
    If $T_\beta$ is the one-root-cut tree, its cut is an internal cut of $T$, as it is followed by $\PureBB{\targetbd' - \beta - \cbd}$, whose root is a branch node as $\targetbd' - \beta - \cbd > \rbd - \beta \ge 0$.
    If $T_\beta$ is the minimal leaf-cut tree, then, its time is at most that of the one-root-cut tree $T'_\beta$ proving $\targetbd' - \beta$, so $\alpha \ge \alpha_\beta$.
    Hence, $\alpha = \alpha_\beta$ and $T'_\beta$ has the same tree time as $T_\beta$.
    Replacing $T_\beta$ by $T'_\beta$ yields a minimal-time tree in which, as before, the cut at the child of $u$ is internal.
    Thus, for $\w(1) = \alpha$, some minimal-time tree proving $\targetbd$ has an internal cut.
\end{proof}

    In \Cref{thm:alpha}, $\targetbd'$ is the remaining bound at the node reached from the root by $k_1$ left and $k_2$ right branches, and $\alpha_\lbd$ and $\alpha_\rbd$ are the values of $\w(1)$ at which the one-root-cut and minimal leaf-cut trees for its two children have the same tree time.
    The condition says that, at $\w(1) = \alpha(\targetbd)$, the one-root-cut tree has time that is no worse than that of the minimal leaf-cut tree for one of these children.
    The forward direction also gives a lower bound: if $u$ is a deepest branch-then-cut node of a minimal-time tree with respect to $\w(1) = \alpha(\targetbd)$, then its child that is not a cut node has a subtree with a cut (by \Cref{lem:impossiblesub}) but no internal cut, so, by the proof of \Cref{lem:pattern}, this subtree is the minimal leaf-cut tree for its target bound, and the same comparison shows that $\alpha(\targetbd)$ is at least the corresponding threshold.
    Hence, for such $k_1$ and $k_2$, we have $\min\{\alpha_\lbd,\alpha_\rbd\} \le \alpha(\targetbd) \le \max\{\alpha_\lbd,\alpha_\rbd\}$.

\begin{figure}
    \centering
    \begin{minipage}[b]{0.24\textwidth}
        \centering
        \scalebox{0.7}{%
        \begin{forest}
            for tree = {big node}
            [0, big cut
                [3
                    [6 [9] [10]]
                    [7 [10] [11]]
                ]
            ]
        \end{forest}%
        }\\[4pt]
        {\small $\treetime(T_1) = 1 + 7a$}
    \end{minipage}%
    \hfill
    \begin{minipage}[b]{0.24\textwidth}
        \centering
        \scalebox{0.7}{%
        \begin{forest}
            for tree = {big node}
            [0
                [3, big cut
                    [6 [9] [10]]
                ]
                [4
                    [7, big cut [10]]
                    [8]
                ]
            ]
        \end{forest}%
        }\\[4pt]
        {\small $\treetime(T_2) = 5 + 4a$}
    \end{minipage}%
    \hfill
    \begin{minipage}[b]{0.24\textwidth}
        \centering
        \scalebox{0.7}{%
        \begin{forest}
            for tree = {big node}
            [0
                [3, big cut
                    [6, big cut [9]]
                ]
                [4
                    [7, big cut [10]]
                    [8]
                ]
            ]
        \end{forest}%
        }\\[4pt]
        {\small $\treetime(T_3) = 5 + 2a + 14/5$}
    \end{minipage}%
    \hfill
    \begin{minipage}[b]{0.24\textwidth}
        \centering
        \scalebox{0.7}{%
        \begin{forest}
            for tree = {big node}
            [0, big cut
                [3, big cut
                    [6 [9] [10]]
                ]
            ]
        \end{forest}%
        }\\[4pt]
        {\small $\treetime(T_4) = 1 + a + 42/5$}
    \end{minipage}
    \caption{Four trees from \Cref{ex:twocut-nonmonotone} proving target bound $\targetbd = 8$ under the setting $(\cbd,\lbd,\rbd) = (3,3,4)$ and $\w = (1, a, 14/5, \infty, \ldots)$. The tree time is shown below each tree, where $14/5 = \w(2)$ for $T_3$, and $42/5 = 3 \w(2)$ for $T_4$.}
    \label{fig:twocut-nonmonotone}
\end{figure}
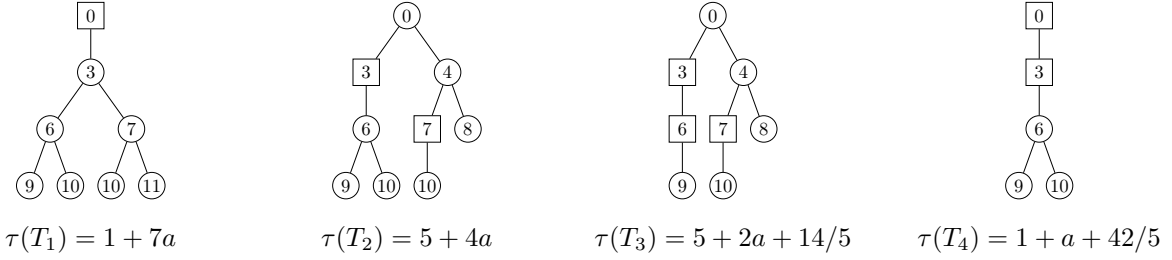

Lastly, we end this subsection with an example showing that similar monotonicity of structure does not follow when more than one cut per root-to-leaf path is allowed.
\begin{example}\label{ex:twocut-nonmonotone}
    Consider a single constant cut family with $\cbd = \lbd = 3$, $\rbd = 4$, target bound $\targetbd = 8$, and $\w = (1, a, 14/5, \infty, \ldots)$ with $a \in [1, 14/5]$, in which at most two cuts are used on any root-to-leaf path.
    \Cref{fig:twocut-nonmonotone} depicts four trees proving $\targetbd$.
    Tree $T_1$ has one root cut.
    Tree $T_2$ has an internal cut in the left subtree and a leaf cut in the right subtree.
    Tree $T_3$ has three leaf cuts.
    Tree $T_4$ has two root cuts.
    As $a = \w(1)$ increases, the unique minimal-time tree goes from $T_1$ to $T_2$ at $a = 4/3$, then $T_3$ at $a = 7/5$, and finally $T_4$ at $a = 8/5$.
    Thus, in contrast to the one-cut setting, cuts first move from the root toward the leaves and then back to the root as $\w(1)$ increases.
\end{example}

A curious consequence of \cref{ex:twocut-nonmonotone} is that increasing $\w(1)$ can lead to a minimal-time tree having strictly more root cuts.
The underlying reason is that, in our model, cuts affect the time of all nodes below them, so increasing $\w(1)$ penalizes having many nodes with exactly one cut above them.
The possible remedies are to either move cuts down the tree or to add cut nodes (to reduce overall tree size), and the latter becomes increasingly attractive as $\w(1)$ approaches $\w(2)$.

\subsection{Effect of Cut Generation Time}
\label{sec:separation}

\Cref{thm:rootcuts,thm:sparse-then-dense} hold for generalized time functions, so cut generation time does not affect the sufficiency of root cuts under symmetric branching or the best way to order consecutive cut nodes.
However, it may affect the number of cut nodes in a minimal-time tree and, under asymmetric branching, where they are placed.

\Cref{prop:separation-monotone} proves the intuitive fact that, as a cut family $i \in [\ncutfams]$ with constant separation time becomes increasingly costly to generate,
minimal-time trees proving the same target bound use weakly fewer cut nodes from that family,
and the remaining tree time, excluding the time to generate these cuts, weakly increases.
For any tree $T$ and family $i \in [\ncutfams]$,
let $n_i(T)$ denote the number of cut nodes from family $i$ in $T$,
and define the tree time without the effect of generating cuts from family $i$ as  
    $\treetime_{-i}(T) \defeq \treetime_{\W}(T) - \sigma_i n_i(T)$,
where $\sigma_i$ is the separation time of family $i$ under $\W$.
Since each cut node of family $i$ contributes exactly $\sigma_i$ to $\treetime_{\W}(T)$, $\treetime_{-i}(T)$ does not depend on $\sigma_i$; that is, $\treetime_{\W}(T) = \treetime_{-i}(T) + \sigma_i n_i(T)$ for every value of $\sigma_i$.

\begin{proposition}
\label{prop:separation-monotone}
    Fix $i \in [\ncutfams]$ as a cut family with constant separation time.
    For $\sigma \ge 0$, let $\W_\sigma$ denote a generalized time function in which family $i$ has separation time $\sigma$, with all other components fixed.
    Let $\sigma_i < \sigma'_i$, and let $T$ and $T'$ be minimal-time trees proving $\targetbd$ with respect to $\W_{\sigma_i}$ and $\W_{\sigma'_i}$, respectively.
    Then $n_i(T') \le n_i(T)$ and $\treetime_{-i}(T') \ge \treetime_{-i}(T)$.
    Moreover, if $n_i(T') = n_i(T)$, then $\treetime_{-i}(T') = \treetime_{-i}(T)$, and both $T$ and $T'$ are minimal-time trees with respect to $\W_{\sigma}$ for every $\sigma \in [\sigma_i, \sigma'_i]$.
\end{proposition}

\begin{proof}
     By minimality, %
     $\treetime_{\W_{\sigma_i}}\!(T) \le \treetime_{\W_{\sigma_i}}\!(T')$
     and
     $\treetime_{\W_{\smash{\sigma'_i}}}\!(T') \le \treetime_{\W_{\smash{\sigma'_i}}}\!(T)$, i.e.,
    \begin{equation*}
        \begin{aligned}
            &\treetime_{-i}(T) + \sigma_i n_i(T) \le \treetime_{-i}(T') + \sigma_i n_i(T')
        \quad \text{and} \\
            &\treetime_{-i}(T') + \sigma'_i n_i(T') \le \treetime_{-i}(T) + \sigma'_i n_i(T).
        \end{aligned}
    \end{equation*}
    Adding these yields $(\sigma'_i - \sigma_i)\big(n_i(T') - n_i(T)\big) \le 0$, so $n_i(T') \le n_i(T)$, and the first inequality then gives $\treetime_{-i}(T) - \treetime_{-i}(T') \le \sigma_i\big(n_i(T') - n_i(T)\big) \le 0$.

    Now suppose that $n_i(T) = n_i(T') = n$. The two inequalities become $\treetime_{-i}(T) \le \treetime_{-i}(T')$ and $\treetime_{-i}(T') \le \treetime_{-i}(T)$, so $\treetime_{-i}(T) = \treetime_{-i}(T')$, and hence $\treetime_{\W_\sigma}(T) = \treetime_{\W_\sigma}(T')$ for all $\sigma$.
    As $T$ is minimal-time for $\sigma_i$ and $T'$ for $\sigma'_i$, 
    both are minimal-time for $\sigma_i$ and for $\sigma'_i$.
    Finally, let $\sigma \in (\sigma_i, \sigma'_i)$ and let $S$ be a minimal-time tree for $\sigma$.
    Applying the first part of the proposition to the tree pairs $(T, S)$ and $(S, T')$ gives $n \ge n_i(S) \ge n$, so $n_i(S) = n$.
    The argument above applied to $(T, S)$ gives $\treetime_{-i}(T) = \treetime_{-i}(S)$.
    Hence, $\treetime_{\W_\sigma}(T) = \treetime_{\W_\sigma}(S) = \treetime_{\W_\sigma}(T')$, so $T$ and $T'$ are minimal-time for $\sigma$.
\end{proof}

  We next consider where cuts are placed when an increase in separation time strictly decreases the number of cut nodes.
  One might expect cuts to be pushed towards the root,
  but this does not hold in general:
  using fewer cut nodes can change where cuts are placed in either direction.
A single root cut may replace several local cuts (\Cref{ex:separation-placement}), but root cuts may also be replaced by fewer, deeper cuts (\Cref{ex:separation-away}).
We write $\w(k) = M$ for a value $M$ large enough that using $k$ or more cuts is never optimal.

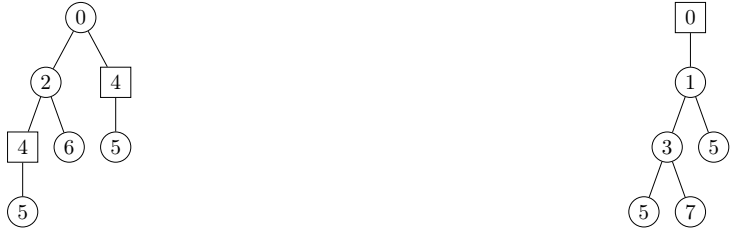
\begin{figure}
    \centering
    \begin{minipage}[t]{0.45\textwidth}
        \centering
        \scalebox{0.8}{%
        \begin{forest}
            for tree = {big node}
            [0
                [2
                    [4, big cut [5]]
                    [6]
                ]
                [4, big cut [5]]
            ]
        \end{forest}%
        }
    \end{minipage}
    \hfil
    \begin{minipage}[t]{0.45\textwidth}
        \centering
        \scalebox{0.8}{%
        \begin{forest}
            for tree = {big node}
            [0, big cut
                [1
                    [3 [5] [7]]
                    [5]
                ]
            ]
        \end{forest}%
        }
    \end{minipage}
    \caption{The trees $T_0$ (left) and $T_1$ (right) from \Cref{ex:separation-placement}. With separation time $\sigma$ per cut node, the minimal-time tree switches from $T_0$ to $T_1$ at $\sigma = 1/5$.}
    \label{fig:separation-placement}
\end{figure}

\begin{example}\label{ex:separation-placement}
    Let $(\cbd, \lbd, \rbd) = (1, 2, 4)$, $\targetbd = 5$, $\w = (1, 1.4, M, \ldots)$,
    for a single constant-strength cut family with constant separation time $\septime_1 \equiv \sigma$.
    Tree $T_0$ in the left panel of \Cref{fig:separation-placement} has two leaf cuts and tree time $5 + 2\w(1) + 2\sigma = 7.8 + 2\sigma$, 
    while tree $T_1$ in the right panel has one root cut and tree time $1 + 5\w(1) + \sigma = 8 + \sigma$.
    The unique minimal-time tree is $T_0$ when $\sigma < 1/5$ and $T_1$ when $1/5 < \sigma < 1$.
    For $\sigma > 1$, the pure branching tree ($T(5)$, with 9 nodes) is optimal.
\end{example}

\begin{example}\label{ex:separation-away}
    \looseness=-1
    Let $(\cbd, \lbd, \rbd) = (3, 3, 7)$, $\targetbd = 12$, $\w = (1, 1.4, 1.9, 2, 3.4, M, \ldots)$, 
    for a single constant-strength cut family with constant separation time $\septime_1 \equiv \sigma$.
    The pure cutting tree, with four cuts, has time $9.7 + 4\sigma$.
    The tree with two root cuts followed by a branch node, whose left child is a leaf cut and right child is a leaf, has time $10.1 + 3\sigma$.
    The unique minimal-time tree is the former when $\sigma < 2/5$, the latter when $2/5 < \sigma < 49/30$, and pure branching for larger $\sigma$, so increasing separation time replaces two root cuts by a single leaf cut.
\end{example}

\Cref{prop:separation-monotone} requires family $i$ to have constant separation time.
When $\septime_i$ depends on the number of prior rounds, two trees with the same cut counts can have different separation times, %
and, in contrast to the second part of \Cref{prop:separation-monotone}, the minimal-time tree can change in structure even when the number of cut nodes does not.

\subsection{Empirical Evaluation of the Root-Cut Approximation}
\label{sec:simulation}

In this section, we empirically assess
(via simulations in our abstract model, using the one-cut setting of \Cref{sec:1cut})
the strength of a root-cut tree and the tightness of the approximation bound $\frac{\tau^*_{\mathrm{root}}(\targetbd)}{\tau^*_\w(\targetbd)}$ in \Cref{prop:notbad}, where $\tau^*_{\mathrm{root}}(\targetbd)$ represents the tree time of the best root-cut tree and $\tau^*_\w(\targetbd)$ represents the tree time of a minimal-time tree.
In the one-cut setting, the best root-cut tree has either no cut or a single root cut, so $\tau^*_{\mathrm{root}}(\targetbd) = \min\{p(\targetbd),\, 1 + \w(1)\,p(\targetbd - \cbd)\}$, and we compute $\tau^*_\w(\targetbd)$ exactly with the recursion in \Cref{apdx:recursion}.
We evaluate $\w(1) \in [1,2]$ in increments of $1/50$, as well as the value of $\w(1)$ at which the ratio peaks, discussed below.%
\footnote{Tree times are computed in exact integer arithmetic when they fit in 64 bits, and in double precision otherwise. Precision (not time) is what limits computation at large depths.}

\looseness=-1
We test settings with $\targetbd = k_1\lbd + k_2\rbd + \cbd$ for some $k_1, k_2 \in \nonnegints$ with $k_1 + k_2 \ge 1$, so the condition of \Cref{lem:existcuts} holds.
We refer to $\lceil \targetbd / \lbd \rceil$, the depth of $\PureBB{\targetbd}$, as the depth of a setting.
We consider three groups of settings, comprising over 55 million parameter combinations:
(A) all $\cbd \le \lbd \le \rbd$ in $[50]$ with $k_1, k_2 \in [0,50]$;
(B) the nearly symmetric settings $(\cbd, 200, \rbd)$ with $\rbd \in \{201, 202\}$, $\cbd \in [100]$, and $k_1, k_2 \in [0,400]$;
and (C) for each $d \in \mathcal{D} \defeq \{5, 10, 50, 100, 500, 1000\}$, all settings $(\cbd, \lbd, \rbd)$ with $\lbd \in \{25, 50, 100, 200\}$, $\rbd - \lbd \in \{1,2\}$, $\cbd \le \min\{\lbd, 100\}$, and depth exactly $d$.
We focus on nearly symmetric settings in (B) and (C) because, in (A), the largest ratio is attained with $\rbd = \lbd + 1$ and does not increase with $\rbd - \lbd$, while the ratio is always $1$ when $\lbd = \rbd$.

\begin{figure}
    \centering
    \ifspringer
        \includegraphics[width=\textwidth]{figs/depth_examples.pdf}
    \else
        \includegraphics[width=\textwidth]{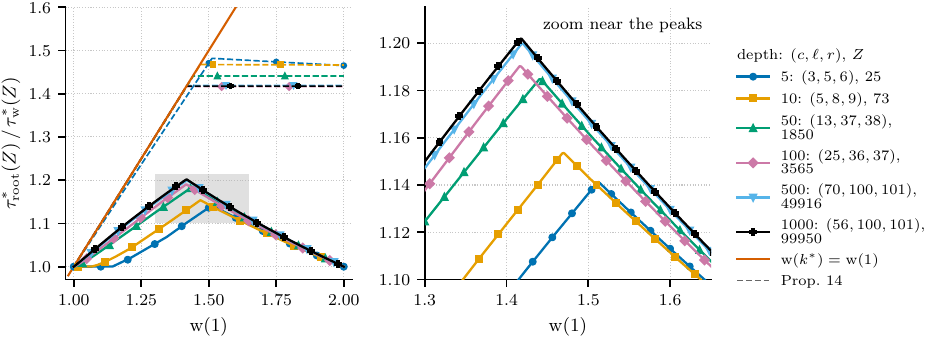}
    \fi
    \caption{Setting at each depth with worst-case approximation ratio as $\w(1)$ varies. The diagonal line is the bound $\w(k^*)$, and the dashed lines represent the bound of \Cref{prop:notbad}.}
    \label{fig:apx}
\end{figure}

\begin{figure}
    \centering
    \includegraphics[width=.8\textwidth]{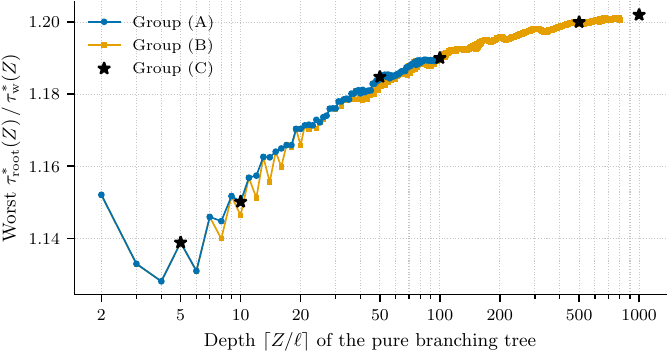}
    \caption{Largest approximation ratio at each depth $\lceil \targetbd/\lbd \rceil$. Groups (A) and (B) are shown up to depth $101$ and $801$, beyond which only settings with $\rbd \gg \lbd$ reach that depth.}
    \label{fig:depthworst}
\end{figure}

\looseness=-1
\Cref{fig:apx} shows, for each $d \in \mathcal{D}$, 
the setting of depth $d$ with the worst (largest) approximation ratio among all settings tested.
Following \Cref{lem:pattern}, the approximation ratio exhibits a stable--increase--decrease pattern as $\w(1)$ gradually increases, maximized at $\w(1) = \frac{p(\targetbd)-1}{p(\targetbd-\cbd)}$,
which is when the pure branching tree and the one-root-cut tree have the same tree time.
\Cref{fig:depthworst} shows how the worst-case ratio grows with depth, from $1.141$ at depth $5$ to $1.191$ at depth $100$ and $1.202$ at depth $1000$.
The bound of \Cref{prop:notbad} provides a conservative estimate of the actual ratio: at the peak, it is between $1.41$ and $1.49$ for the settings in \Cref{fig:apx}.
Hence, a root-cut tree may be a practically-effective (and simpler, computationally) alternative to a minimal-time tree in our model.

\begin{figure}
    \setlength{\subplotwidth}{0.45\textwidth}  %
    \centering
    \begin{minipage}{\subplotwidth}
        \centering
        \includegraphics[width=\linewidth]{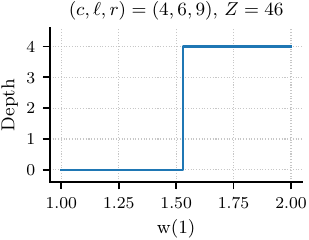}
    \end{minipage}
    \hfill
    \begin{minipage}{\subplotwidth}
        \centering
        \includegraphics[width=\linewidth]{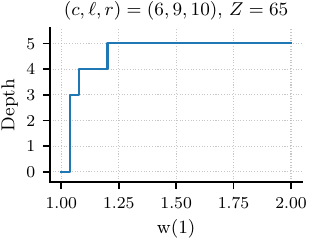}
    \end{minipage}

    \vspace{0.3cm}

    \begin{minipage}{\subplotwidth}
        \centering
        \includegraphics[width=\linewidth]{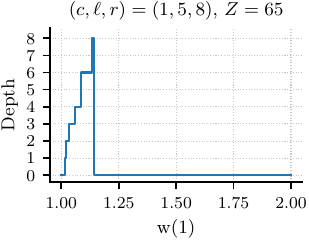}
    \end{minipage}
    \hfill
    \begin{minipage}{\subplotwidth}
        \centering
        \includegraphics[width=\linewidth]{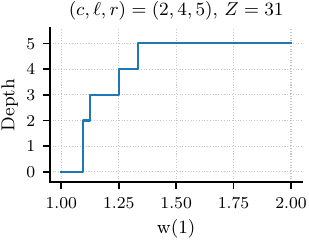}
    \end{minipage}
    \caption{Depth of shallowest internal cut over all minimal-time trees as $\w(1)$ varies. A value of $0$ indicates that no minimal-time tree has an internal cut.}
    \label{fig:subplots}
\end{figure}

Lastly, as $\w(1)$ grows, by \Cref{lem:pattern}, the minimal-time structure shifts from root-cut to internal-cut to leaf-cut.
Intuitively, cuts get pushed down in a minimal-time tree.
\emph{Does this shift progress one layer at a time?}
The answer is \emph{not always}.
Empirically, the shallowest internal cut among all minimal-time trees may jump from the root to a much deeper layer as $\w(1)$ grows, as depicted in \Cref{fig:subplots}.
Specifically, these plots show the depth of the shallowest internal cut over all minimal-time trees, for four settings and incrementing $\w(1)$ by thousandths.
A value of $0$ indicates that no minimal-time tree has an internal cut, i.e., all cuts are at the root or at the leaves.

\section{Conclusions}
\label{sec:conclusions}

\looseness=-1
This paper highlights the structural importance of root cuts within our abstract \bc framework.
In particular, we identify settings, both symmetric and asymmetric, in which some minimal-time tree has no internal cuts, even when accounting for cut generation time.
More generally, we bound how much the best root-cut tree can deviate from a minimal-time tree, and we prove that, for constant cuts, the best root-cut tree is asymptotically optimal.
In addition, we characterize when internal cuts are necessary in a one-cut setting, and we show that cut generation time can move cuts either toward or away from the root.
 
While this model captures a phenomenon commonly observed in state-of-the-art solvers, that cuts are used aggressively at the root, it intentionally abstracts away complex solver interactions.
Specifically, cutting and branching interact solely through bound improvement and the slowdown in node processing time, leaving out other mechanisms such as pruning, restarts, primal heuristics, and the sequence of fractional solutions that a solver encounters.
Consequently, preliminary computational experiments evaluating cut placement only exhibit noise-level differences in performance.
This motivates generalizations of the current model.
For example, as branching decisions can be significantly affected by cuts in practice, a natural extension is to incorporate cut-dependent branching decisions.

\ifllncs
  \begin{credits}
  \subsubsection{\ackname}
  This material is based upon work supported by the Air Force Office of Scientific Research under award number FA9550-23-1-0340.

  \subsubsection{\discintname}
  The authors have no competing interests to declare. %
  \end{credits}
\else
    \ifspringer
         \begin{acknowledgements}
            We appreciate the constructive feedback from IPCO reviewers that helped improve our paper.
            Based on initial results by the authors, Claude's Fable 5.1 and Opus 5.5 models were used to check and improve proofs, correct grammar, fix notational inconsistencies, and create better figures.
            Afterwards, these models helped simplify some results and generalize others, including an initial version of the asymptotic analysis at the authors' prompting, to merge what was previously two sections into one (\cref{sec:rootcut-quality}) and to precisely consider the effect of separation time (\cref{sec:separation}).
            Finally, the same models were used to generate the code for the simulation experiments and reproducible figures, replacing our previous slower mixed-integer optimization approach.
         \end{acknowledgements}
        
        \begingroup
        \renewcommand{\ackname}{Declarations}
        \begin{acknowledgements}
            \emph{Competing interests and funding.} The authors have no competing interests to declare. This material is based upon work supported by the Air Force Office of Scientific Research under award number FA9550-23-1-0340.
            \emph{Code availability.} The code for the simulations in \cref{sec:simulation} is available at \url{https://github.com/akazachk/abstract-branch-and-cut}.
        \end{acknowledgements}
        \endgroup
    \else
      \paragraph{Acknowledgments.}
      We appreciate the constructive feedback from IPCO reviewers that helped improve our paper.
      Based on initial results by the authors, Claude's Fable 5.1 and Opus 5.5 models were used to check and improve proofs, correct grammar, fix notational inconsistencies, and create better figures.
      Afterwards, these models helped simplify some results and generalize others, including an initial version of the asymptotic analysis at the authors' prompting, to merge what was previously two sections into one (\cref{sec:rootcut-quality}) and to precisely consider the effect of separation time (\cref{sec:separation}).
      Finally, the same models were used to generate the code for the simulation experiments and reproducible figures, replacing our previous slower mixed-integer optimization approach.
    
      \paragraph{Funding.} 
      This material is based upon work supported by the Air Force Office of Scientific Research under award number FA9550-23-1-0340.
    
      \paragraph{Code availability.} The code for the simulations in \cref{sec:simulation} is available at \url{https://github.com/akazachk/abstract-branch-and-cut}.
    \fi
\fi

\ifllncs
  \bibliographystyle{splncs04}
\else
  \bibliographystyle{plainnat}
\fi
\bibliography{include/akazachk}

\appendix 

\section*{Appendix}
\addcontentsline{toc}{section}{Appendix}

\setcounter{section}{0}
\renewcommand\thesection{\Alph{section}}

\section{Exact Approaches to Computing Minimal-Time One-Cut Trees}

In this section, we show how to compute a minimal-time one-cut tree in two ways, first with a recursive formula, and then with a mixed-integer linear program.

\subsection{Exact Recursion}
\label{apdx:recursion}

In the one-cut setting, the minimal tree time can be computed by a recursion on the remaining bound.
Consider a node with remaining bound $z$ and no cut above it.
If it is a cut node, then no further cut can be added below it, so its subtree is the cut node followed by $\PureBB{z - \cbd}$, with tree time $1 + \w(1)\,p(z - \cbd)$.
If it is a branch node, then, by \cite[Lemma~24]{KazLeBSan24_abstractbc-mpb}, its two subtrees are minimal-time trees for remaining bounds $z - \lbd$ and $z - \rbd$.
Hence, the minimal tree time $\tau^*_\w(z)$ satisfies
\begin{equation*}
    \tau^*_\w(z) =
    \begin{cases}
        1 & \text{if } z \le 0, \\
        1 + \min\bigl\{\tau^*_\w(z - \lbd) + \tau^*_\w(z - \rbd),\; \w(1)\,p(z - \cbd)\bigr\} & \text{otherwise},
    \end{cases}
\end{equation*}
where $p(z) = 1$ for $z \le 0$ and $p(z) = 1 + p(z - \lbd) + p(z - \rbd)$ otherwise.
The best root-cut tree has either no cut or a single root cut, so its tree time is $\min\{p(\targetbd),\, 1 + \w(1)\,p(\targetbd - \cbd)\}$, and the approximation ratio in \Cref{fig:apx} is this value divided by $\tau^*_\w(\targetbd)$.
Recording which term attains the minimum at each $z$ recovers a minimal-time tree, and keeping all terms that attain it recovers all minimal-time trees; for example, the depth of the shallowest internal cut over all minimal-time trees, as in \Cref{fig:subplots}, is computed by the same recursion.

Every remaining bound that appears in the recursion has the form $\targetbd - k_1\lbd - k_2\rbd$, so, with memoization, $\tau^*_\w$ and $p$ are evaluated at most $\targetbd$ times each at positive arguments, and the recursion takes $O(\targetbd)$ arithmetic operations for each value of $\w(1)$.
Moreover, $\tau^*_\w(z)$ does not depend on $\targetbd$, so a single pass of the recursion for a fixed $(\cbd,\lbd,\rbd)$ gives the minimal tree time for every target bound of that triple at once.
Writing $\w(1) = k/50$ for an integer $k$, all tree times scaled by $50$ are integers, which we compute exactly whenever they fit in 64 bits, so that ties between minimal-time trees, such as at $\w(1) = \alpha(\targetbd)$ or $\w(1) = 2$, are detected exactly rather than up to a solver tolerance.
As $p$ grows exponentially (for example, $p(1000)$ has 210 digits when $(\lbd,\rbd) = (1,2)$), larger values are computed in double precision, which was needed for about $73\%$ of the settings in \Cref{sec:simulation}, all of them in groups (A) and (B).
The reported worst-case settings and their peak ratios are rechecked with exact integer and rational arithmetic.
Evaluating all values $\w(1) \in \{1.00, 1.02, \ldots, 2.00\}$ at once, the full set of $55{,}359{,}582$ settings in \Cref{sec:simulation} took $224$ seconds with a Python implementation using $10$ parallel workers on an Apple M1 Ultra processor.

\subsection{Mixed-Integer Linear Programming Formulation}\label{apdx:formulation}

Another way to compute the trees in \cref{sec:simulation} is by solving a mixed-integer linear program over a complete binary tree of candidate nodes.
The input is a one-cut setting $(\targetbd, \cbd, \lbd, \rbd)$ and the value $\w(1)$.
Let $\nodes$ be the node set of a complete binary tree with root $0$ and depth $\lceil \targetbd/\lbd \rceil$, which is the maximum necessary depth when considering only branch nodes.
For $v \in \nodes \setminus \{0\}$, let $\mathrm{par}(v)$ denote the parent of $v$, and let $\nodes_L$ and $\nodes_R$ be the sets of left and right children in $\nodes$.
For $v \in T$, let $P(v)$ denote the nodes on the path from the root to $v$, including $v$.

For each $v \in \nodes$, define a binary variable $y_v$ to indicate whether $v$ is a node of the tree, and a binary variable $x_v$ for whether a cut node is placed immediately above $v$, i.e., between $v$ and its parent, or above the root if $v = 0$.
Let $z_v \in \nonnegreals$ be the variable for the bound at $v$.
The binary variable $q_v$ indicates whether a cut is applied on the path from the root to $v$, in which case $v$ has time $\w(1)$, and the binary variable $m_v$ linearizes the product $q_v y_v$.
The formulation is
\begin{subequations}\label{eq:milp}
\begin{align}
    \min_{x,y,q,m,z} \quad & \sum_{v \in \nodes} \big(x_v + y_v + (\w(1) - 1)\,m_v\big) \nonumber \\
    \quad
    & z_v = z_{\mathrm{par}(v)} + \lbd + \cbd\, x_v && v \in \nodes_L, \label{eq:milp-left}\\
    & z_v = z_{\mathrm{par}(v)} + \rbd + \cbd\, x_v && v \in \nodes_R, \label{eq:milp-right}\\
    & z_0 = \cbd\, x_0, \quad y_0 = 1, \label{eq:milp-root}\\
    & y_v \ge x_v && v \in \nodes, \label{eq:milp-cutexists}\\
    & y_v = y_{v'} && \text{siblings } v, v' \in \nodes, \label{eq:milp-siblings}\\
    & y_{\mathrm{par}(v)} \ge y_v \ge 1 - z_{\mathrm{par}(v)}/\targetbd && v \in \nodes \setminus \{0\}, \label{eq:milp-exists}\\
    & \textstyle\sum_{u \in P(v)} x_u \le 1 && v \in \nodes, \label{eq:milp-onecut}\\
    & q_v \ge \textstyle\sum_{u \in P(v)} x_u && v \in \nodes, \label{eq:milp-q}\\
    & m_v \le q_v, \quad m_v \le y_v, \quad m_v \ge q_v + y_v - 1 && v \in \nodes, \label{eq:milp-mccormick}\\
    & x_v, y_v, q_v, m_v \in \{0,1\}, \quad z_v \in \nonnegreals && v \in \nodes. \nonumber
\end{align}
\end{subequations}

The objective is the tree time.
Each node of the tree contributes $1$, plus $\w(1) - 1$ if a cut lies on its path from the root, and each cut node contributes $1$, as it has no other cut above it in the one-cut setting.
Constraints \eqref{eq:milp-left}--\eqref{eq:milp-root} compute the bound at each node, where a cut above $v$ adds $\cbd$.
By \eqref{eq:milp-cutexists}, a cut is only placed above a node of the tree, and by \eqref{eq:milp-siblings}, either both children of a node are in the tree or neither is.
By \eqref{eq:milp-exists}, a node is in the tree only if its parent is, and the children of a node must be in the tree whenever its bound is less than $\targetbd$, so every leaf proves $\targetbd$.
Constraint \eqref{eq:milp-onecut} allows at most one cut on each root-to-leaf path, and \eqref{eq:milp-q} sets $q_v = 1$ when there is such a cut above $v$.
Finally, \eqref{eq:milp-mccormick} are the McCormick inequalities for $m_v = q_v y_v$; as $\w(1) > 1$, an optimal solution has $m_v = 1$ exactly when $v$ is in the tree and has a cut above it.
For the best root-cut tree, we also fix $x_v = 0$ for all $v \in \nodes \setminus \{0\}$.

\end{document}